\documentclass[a4paper,11pt]{article}
\usepackage{amsmath,amsthm,amssymb,dsfont,graphicx,xspace,epsfig}
\usepackage[dvipsnames]{xcolor}
\usepackage[T1]{fontenc}
\usepackage[plain]{fullpage}
\usepackage{thm-restate}
\usepackage[hidelinks]{hyperref}
\usepackage{bm}
\usepackage{color}
\usepackage{comment}
\usepackage{mathrsfs}
\usepackage{cleveref}
\usepackage[shortlabels]{enumitem}
\usepackage{float}
\usepackage{dutchcal}
\usepackage{pict2e}
\usepackage{mathtools}
\usepackage{tikz}
\usepackage{subfigure}
\usepackage{xstring,refcount}
\usepackage{authblk}
\usepackage{yfonts}

\usetikzlibrary{fit, backgrounds, matrix, arrows.meta}
\usetikzlibrary{calc,arrows,positioning}
\usetikzlibrary{decorations}
\usetikzlibrary{decorations.pathmorphing}
\usetikzlibrary{decorations.pathreplacing}
\usetikzlibrary{decorations.shapes}
\usetikzlibrary{decorations.text}
\usetikzlibrary{decorations.markings}
\usetikzlibrary{decorations.fractals}
\usetikzlibrary{decorations.footprints}

\tikzset{vertex/.style = {circle,fill=black,minimum size=6pt, inner sep=0pt, outer sep=3pt}}
\tikzset{svertex/.style = {circle,fill=black,minimum size=5pt, inner sep=0pt, outer sep=2pt}}
\tikzset{labelledvertex/.style = {circle,fill=none,draw,very thick, inner sep=2pt, outer sep=2pt}}
\tikzset{Xrect/.style = {rectangle,fill=none,draw,very thick, minimum width=1.2cm, minimum height=0.8cm, outer sep=0pt}}
\tikzset{edge/.style = {thick,dashed,->,> = latex,dash pattern={on 5pt off 4pt}}}
\tikzset{arc/.style = {thick,->,> = latex}}
\tikzset{sarc/.style = {thin,->,> = latex}}
\tikzset{vsarc/.style = {very thin,->,> = latex}}
\tikzset{digon/.style = {thick,<->,> = latex}}
\tikzset{bigvertex/.style = {shape=circle,draw}}

\definecolor{g-blue}{rgb}{0.0, 0.5, 1.0}
\definecolor{g-green}{rgb}{0.4, 0.9, 0.4522}

\newtheorem{theorem}{Theorem}
\newtheorem{lemma}[theorem]{Lemma}
\newtheorem{corollary}[theorem]{Corollary}
\newtheorem{proposition}[theorem]{Proposition}
\newtheorem{conjecture}[theorem]{Conjecture}

\newtheorem{claim}{Claim}[theorem]
\newcounter{oldthm}{}
\def\@extractthmnum#1.#2{#2}
\newenvironment{claimlabeled}[1][\@nil]{%
    \setcounter{oldthm}{\value{theorem}}%
    \def\tmp{#1}%
    \ifx\tmp\@nnil
    \else
       ~\refused{#1}
        \setcounter{theorem}{\getrefnumber{#1}}
    \fi
    \begin{claim}%
}{%
    \ifx\tmp\@nnil
    \else
        \setcounter{theorem}{\value{oldthm}}
    \fi
    \end{claim}
}

\theoremstyle{definition}

\newtheorem{problem}[theorem]{Problem}

\newenvironment{proofclaim}[1][]{\par\noindent {\it Proof of claim}. }{ \hfill$\lozenge$\par\addvspace{6pt plus 6pt}}

\newcommand{\opentriangle}{
  \raisebox{0.2pt}{\makebox[0.77778em]{
    \setlength{\unitlength}{0.6em}
    \linethickness{0.4pt}
    \begin{picture}(1,1)
    \polygon(0,0)(1,0)(1,1)
    \end{picture}
  }}
}

\newcommand{\defproblem}[3]{
 \vspace{3mm}
\noindent\fbox{
 \begin{minipage}{0.96\textwidth}
 \begin{tabular*}{\textwidth}{@{\extracolsep{\fill}}lr} #1  \\ \end{tabular*}
 {\bf{Input:}} #2 \\
 {\bf{Question:}} #3
 \end{minipage}
 }
 \vspace{3mm}
}

\newcommand{\Ra}{\Rightarrow}
\newcommand{\dic}{\vec{\chi}}
\newcommand{\chia}{\chi_{\rm a}}
\newcommand{\dict}{\vec{\chi}_{\rm a}}
\makeatletter
\newcommand{\overleftrightsmallarrow}{\mathpalette{\overarrowsmall@\leftrightarrowfill@}}
\newcommand{\overrightsmallarrow}{\mathpalette{\overarrowsmall@\rightarrowfill@}}
\newcommand{\overleftsmallarrow}{\mathpalette{\overarrowsmall@\leftarrowfill@}}
\newcommand{\overarrowsmall@}[3]{%
  \vbox{%
    \ialign{%
      ##\crcr
      #1{\smaller@style{#2}}\crcr
      \noalign{\nointerlineskip}%
      $\m@th\hfil#2#3\hfil$\crcr
    }%
  }%
}
\def\smaller@style#1{%
  \ifx#1\displaystyle\scriptstyle\else
    \ifx#1\textstyle\scriptstyle\else
      \scriptscriptstyle
    \fi
  \fi
}
\makeatother
\newcommand{\bid}[1]{\overleftrightsmallarrow{#1}}

\newcommand{\PP}{\mathbb{P}}

\newcommand{\Ccal}{\mathcal{C}}

\newcommand{\Ascr}{\mathscr{A}}
\newcommand{\Cscr}{\mathscr{C}}

\newcommand{\Pscr}{\mathscr{P}}

\newcommand{\forb}[1]{\mathfrak{F}(#1)}

\renewcommand{\epsilon}{\varepsilon}
\renewcommand{\emptyset}{\varnothing}

\renewcommand{\phi}{\varphi}
\let\le\leqslant
\let\ge\geqslant
\let\leq\leqslant
\let\geq\geqslant

\hypersetup{
    colorlinks,
    linkcolor={red!50!black},
    citecolor={blue!50!black},
    urlcolor={blue!80!black}
}
\newcommand{\2}{\vspace{2mm}}

\title{Acyclic dichromatic number of oriented graphs}

\author[1,2]{Jørgen Bang-Jensen}
\author[,3]{Lucas Picasarri-Arrieta\thanks{Research supported by JST ASPIRE JPMJAP2302.}}
\author[1,4]{Anders Yeo}

\affil[1]{Department of Mathematics and Computer Science, University of Southern Denmark, Odense, Denmark}
\affil[2]{School of Mathematics, Shandong University, Jinan 250100, China}
\affil[3]{National Institute of Informatics, Tokyo, Japan}
\affil[4]{Department of Mathematics and Applied Mathematics, University of Johannesburg, Auckland Park, 2006 South Africa }

\date{}

\begin{document}

\maketitle

\begin{abstract}
The dichromatic number $\dic(D)$ of a digraph $D=(V,A)$ is the minimum number of sets in a partition $V_1,\ldots{},V_k$ of $V$ into $k$ subsets so that the induced subdigraph $D[V_i]$ is acyclic for each $i\in [k]$. This is a generalization of the chromatic number for undirected graphs as a graph has chromatic number at most $k$ if and only if the complete biorientation $\bid{G}$ of $G$ (replace each edge by a directed 2-cycle) has dichromatic number at most $k$. In this paper we introduce the acyclic dichromatic number $\dict{}(D)$ of a digraph $D$ as the minimum number of sets in a partition $V_1,\ldots{},V_k$ of $V$ into $k$ subsets so that the induced subdigraph $D[V_i]$ is acyclic for each $i\in [k]$ and each of the bipartite induced subdigraphs $D[V_i,V_j]$ is acyclic for each $1\leq i<j\leq k$. This parameter, which resembles the definition of acyclic chromatic number for undirected graphs, has apparently not been studied before. 
We derive a number of results which display the difference between the dichromatic number and the acyclic dichromatic number, in particular, there are digraphs $D$ with arbitrarily large $\dict{}(D)-\dic{}(D)$, even among tournaments with dichromatic number 2 and bipartite tournaments (where the dichromatic number is always 2). We prove several complexity results, including that deciding whether $\dict{}(D)\leq 2$ is NP-complete already for bipartite digraphs, while it is polynomial for tournaments (contrary to the case for dichromatic number).
We also generalize the concept of heroes of a tournament to acyclic heroes of tournaments.

\medskip

\noindent{}{\bf Keywords:} Dichromatic number; acyclic dichromatic number; tournaments; bipartite tournaments; split digraphs; NP-completeness; acyclic heroes.
\end{abstract}

\section{Introduction}

Terminology  not provided below is consistent with \cite{bang2009}.

Given an integer $k\in \mathbb{N}$, a {\bf $k$-colouring} of a digraph $D$ is a mapping $\phi \colon V(D) \to [k]$. The sets $V_1,\dots,V_k$ where $V_i = \phi^{-1}(i)$ are called the {\bf colour classes} of $\phi$.
We say that $\phi$ is a {\bf $k$-dicolouring} if, for every $i\in [k]$, $D[V_i]$ is acyclic. The {\bf dichromatic number} of $D$ is the least integer $k\in \mathbb{N}$ such that $D$ admits an acyclic $k$-dicolouring.
The notions of dicolouring and dichromatic number were introduced by Erd\H{o}s and Neumann-Lara in the late 1970s~\cite{erdosPNCN1979,neumannlaraJCT33}, rediscovered independently by Mohar~\cite{bokalJGT46,moharJGT43} in the 2000s, and gained a lot of attention since then. It turned out that many classical results on graph colouring admit a generalisation to the directed setting involving the dichromatic number, see for instance~\cite{aboulkerEJC28,andresJGT79,harutyunyanArxiv25,harutyunyanSIDMA25,kawarabayashiSODA25}.

In this paper, we introduce the notion of {\bf acyclic dicolouring} of oriented graphs (that is, digraphs with at most one arc between any two vertices). A $k$-dicolouring $\phi$, with colour classes $V_1,\dots,V_k$, is {\bf acyclic} if, for every $1\leq i<j\leq k$, the bipartite oriented graph $D[V_i,V_j]$ is acyclic. 
The {\bf acyclic dichromatic number} of $D$, denoted $\dict(D)$, is the least integer $k\in \mathbb{N}$ such that $D$ admits a acyclic $k$-dicolouring. Note that $\dict(D)$ is well-defined, as every oriented graph $D$ admits a trivial acyclic $|V(D)|$-dicolouring.

This notion is a directed analogue of the one of {\bf acyclic colouring}~\cite{grunbaumIJM14}, which is a colouring of an undirected graph such that any two adjacent vertices receive distinct colours, and the union of any two colour classes induces a forest ({\it i.e.} is acyclic). In particular, if $D$ is an oriented graph with underlying graph $G$, then any acyclic colouring of $G$ yields an acyclic dicolouring of $D$, and in particular $\dict(D) \leq \chia(G)$, where $\chia(G)$ denotes the acyclic chromatic number of $G$.

Adding the restriction that also the bipartite subdigraphs $D[V_i,V_j]$ are acyclic for each $1\leq i<j\leq k$ may force us to use (many) more colours in an acyclic dicolouring of a digraph $D$ than if we just require each $D[V_i]$ to be acyclic. This also has the effect that the complexity of deciding whether $\dict{(D)}\leq k$ may differ from that of deciding whether $\dic(D)\leq k$, namely, the first problem may be polynomial for some $k$ even if the second one is NP-complete. Indeed, as a prominent case, we prove that this is the case for tournaments when $k=2$.

The paper is organized as follows. In Section~\ref{sec:relation}, we show how to construct (bipartite) digraphs with arbitrarily high acyclic dichromatic number, showing that the difference between the acyclic dichromatic number and the dichromatic number can be arbitrarily high even for bipartite digraphs. We use this to show that the same holds also for the class of tournaments of dichromatic number 2. In Section~\ref{sec:acyclichero}, we generalize the concept of heroes of tournaments to what we call acyclic heroes for tournaments, that is, tournaments so that every tournament that does not contain one of these as a subtournament has bounded acyclic dichromatic number. We formulate a conjecture describing which tournaments are exactly the acyclic heroes, and provide several results which support this conjecture. In Section~\ref{sec:complexity}, we first use the construction from Section~\ref{sec:relation} to obtain an easy proof that determining whether $\dict(D)\leq k$ is NP-complete for every $k$, even for bipartite digraphs. We also show that the same conclusion holds for split digraphs. We then show that, contrary (assuming P is not equal to NP) to the case for the dichromatic number, one can check in polynomial time whether a tournament has acyclic dichromatic number at most 2. In Section~\ref{sec:tightbounds}, we derive tight upper bounds for the acyclic dichromatic number of tournaments and in Section~\ref{sec:degenerate} we show that every 2-degenerate oriented graph has acyclic dichromatic number at most 2 while there are 3-degenerate oriented graphs with arbitrarily high acyclic dichromatic number. Finally, in Section~\ref{sec:remarks} we pose a number of open problems on acyclic dichromatic numbers of oriented graphs.

\section{A construction of oriented graphs with large acyclic dichromatic number}
\label{sec:relation}

It is well-known that the random tournament on $n$ vertices has dichromatic number $\Theta(n/\log n)$ with high probability, see {\it e.g.}~\cite[Section~4]{bensmailJGT87}. In Theorem~\ref{thm:lower_bound_T_n} below, we show that the random tournament on $n$ vertices has acyclic dichromatic number $n - \Theta(\log n)$ with high probability, so there exist tournaments $T$ of order $n$ whose gap between  $\dict{(T)}$ and $\dic{(T)}$ is $n(1-o(1))$. 

\medskip

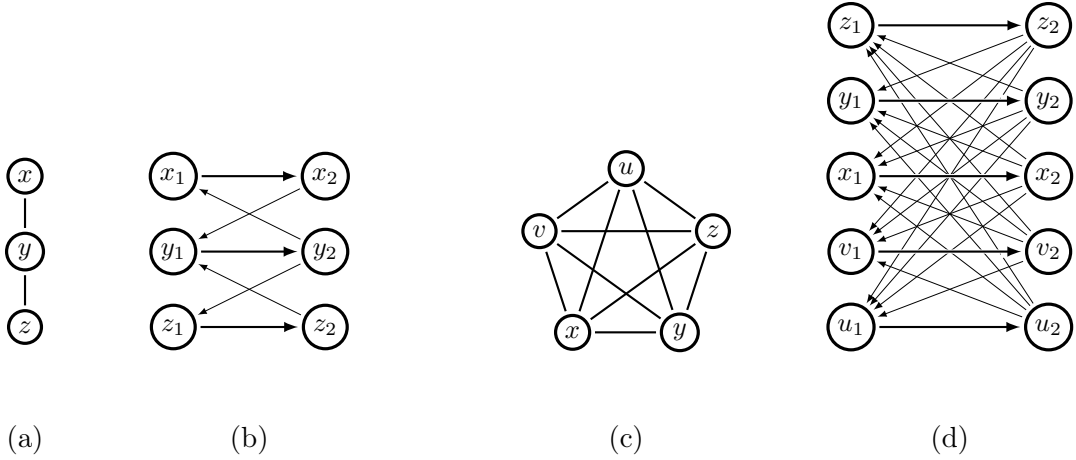
\begin{figure}[ht]
    \centering
    \begin{tikzpicture}
      \node at (0,-1.5) {(a)};
      \node (x) at (0,2) [labelledvertex]{\small $x$};
      \node (y) at (0,1) [labelledvertex]{\small $y$}; 
      \node (z) at (0,0) [labelledvertex]{\small $z$};
      \draw[thick] (x) -- (y);
      \draw[thick] (y) -- (z);
    \end{tikzpicture} \hspace{1cm}
    \begin{tikzpicture}
      \node at (1,-1.5) {(b)};
      \node (x1) at (0,2) [labelledvertex]{\small $x_1$};
      \node (y1) at (0,1) [labelledvertex]{\small $y_1$};
      \node (z1) at (0,0) [labelledvertex]{\small $z_1$};
    
      \node (x2) at (2,2) [labelledvertex]{\small $x_2$};
      \node (y2) at (2,1) [labelledvertex]{\small $y_2$};
      \node (z2) at (2,0) [labelledvertex]{\small $z_2$};
    
      \draw[arc] (x1) to (x2);
      \draw[arc] (y1) to (y2);
      \draw[arc] (z1) to (z2);
    
      \draw[vsarc] (x2) to (y1);
      \draw[vsarc] (y2) to (x1);
      \draw[vsarc] (y2) to (z1);
      \draw[vsarc] (z2) to (y1);
    \end{tikzpicture} \hspace{2cm}
    \begin{tikzpicture}
      \node at (0,-2.4) {(c)};
      \node (u) at (90:1.2) [labelledvertex]{\small$u$};
      \node (v) at (162:1.2) [labelledvertex]{\small$v$};
      \node (x) at (234:1.2) [labelledvertex]{\small$x$};
      \node (y) at (306:1.2) [labelledvertex]{\small$y$};
      \node (z) at (18:1.2) [labelledvertex]{\small$z$};
      \draw[thick] (u) -- (v);
      \draw[thick] (u) -- (x);
      \draw[thick] (u) -- (y);
      \draw[thick] (u) -- (z);
      \draw[thick] (v) -- (x);
      \draw[thick] (v) -- (y);
      \draw[thick] (v) -- (z);
      \draw[thick] (x) -- (y);
      \draw[thick] (x) -- (z);
      \draw[thick] (y) -- (z);
    \end{tikzpicture} \hspace{1cm}
    \begin{tikzpicture}
      \node at (1,-1.5) {(d)};
      \node (u1) at (-0.3,0) [labelledvertex]{\small$u_1$};
      \node (v1) at (-0.3,1) [labelledvertex]{\small$v_1$};
      \node (x1) at (-0.3,2) [labelledvertex]{\small$x_1$};
      \node (y1) at (-0.3,3) [labelledvertex]{\small$y_1$};
      \node (z1) at (-0.3,4) [labelledvertex]{\small$z_1$};
    
      \node (u2) at (2.3,0) [labelledvertex]{\small$u_2$};
      \node (v2) at (2.3,1) [labelledvertex]{\small$v_2$};
      \node (x2) at (2.3,2) [labelledvertex]{\small$x_2$};
      \node (y2) at (2.3,3) [labelledvertex]{\small$y_2$};
      \node (z2) at (2.3,4) [labelledvertex]{\small$z_2$};

      \draw[vsarc] (u2) to (v1);
      \draw[vsarc] (u2) to (x1);
      \draw[vsarc] (u2) to (y1);
      \draw[vsarc] (u2) to (z1);
    
      \draw[vsarc] (v2) to (u1);
      \draw[vsarc] (v2) to (x1);
      \draw[vsarc] (v2) to (y1);
      \draw[vsarc] (v2) to (z1);
    
      \draw[vsarc] (x2) to (u1);
      \draw[vsarc] (x2) to (v1);
      \draw[vsarc] (x2) to (y1);
      \draw[vsarc] (x2) to (z1);
    
      \draw[vsarc] (y2) to (u1);
      \draw[vsarc] (y2) to (v1);
      \draw[vsarc] (y2) to (x1);
      \draw[vsarc] (y2) to (z1);
    
      \draw[vsarc] (z2) to (u1);
      \draw[vsarc] (z2) to (v1);
      \draw[vsarc] (z2) to (x1);
      \draw[vsarc] (z2) to (y1);
      
      \draw[arc, ultra thick, white] (u1) to (u2);
      \draw[arc, ultra thick, white] (v1) to (v2);
      \draw[arc, ultra thick, white] (x1) to (x2);
      \draw[arc, ultra thick, white] (y1) to (y2);
      \draw[arc, ultra thick, white] (z1) to (z2);
      
      \draw[arc] (u1) to (u2);
      \draw[arc] (v1) to (v2);
      \draw[arc] (x1) to (x2);
      \draw[arc] (y1) to (y2);
      \draw[arc] (z1) to (z2);
    \end{tikzpicture}

    \caption{Figures (b) and (d) illustrate the digraph $D(G)$ for the two graphs shown in Figures~(a) and (c), respectively.}\label{fig:vertexsplitG}
\end{figure}

In this section, we provide a construction of tournaments with dichromatic number $2$ but arbitrarily large acyclic dichromatic number that will be used later on.
For any undirected graph $G=(V,E)$, the corresponding symmetric digraph $\bid{G}$ is obtained by replacing each edge of $G$ by a directed 2-cycle. We let $D(G)$ be the oriented graph obtained from $\bid{G}$ by performing the operation called {\bf vertex splitting} in \cite[page 134]{bang2009}, that is, each vertex $v$ is replaced by two vertices $v_{in},v_{out}$ and the arc $v_{in}v_{out}$ and every arc $uv$ is replaced by  the arc $u_{out}v_{in}$. The digraph $D(G)$ is bipartite with bipartition $(V_{in},V_{out})$, where $V_{in}=\{v_{in}:v\in V\}$ and $V_{out}=\{v_{out}:v\in V\}$, and the only arcs from $V_{in}$ to $V_{out}$ is the perfect matching $\{v_{in}v_{out}:v\in V\}$. Observe that every edge $uv$ of $G$ corresponds to the directed 4-cycle $(u_{in},u_{out},v_{in},v_{out})$ and that the only induced cycles of $D(G)$ are 4-cycles of the form above. See Figure~\ref{fig:vertexsplitG} for an example.

It turns out that $\dict(D(G))$ grows with $\chi(G)$.

\begin{proposition}
    \label{prop:chi_DG}
    For every graph $G$, $\dict(D(G)) = \left\lceil \sqrt{\chi(G)}\right\rceil$.
\end{proposition}
\begin{proof}
    Let $k= \dict(D(G))$. 
    We first prove that $k \geq \lceil \sqrt{\chi(G)}\rceil$.
    To see this, let $\phi\colon V(D(G))\rightarrow [k]$ be an acyclic dicolouring of $D(G)$ and assign each vertex $v\in V(G)$ the colour $\psi(v)=(\phi(v_{in}),\phi(v_{out}))$. We claim that $\psi$ is a proper colouring of $G$. To see this, it suffices to see that if $uv\in E$ is an edge of $G$ such that $\psi(u)=\psi(v)$, then either all vertices of the 4-cycle $(u_{in},u_{out},v_{in},v_{out})$  received the same colour $i\in [k]$ by $\phi$ or the vertices $u_{in},v_{in}$ received one colour $i$ under $\phi$ and $u_{out},v_{out}$ received another colour $j$ under $\phi$. In the first case $\phi^{-1}(i)$ is not acyclic and in the second case the 4-cycle $(u_{in},u_{out},v_{in},v_{out})$ alternates between the sets $\phi^{-1}(i)$ and $\phi^{-1}(j)$ and both cases contradict that $\phi$ is an acyclic dicolouring of $D(G)$. Since $\psi$ uses $k^2$ colours, this shows that $k\geq \sqrt{\chi(G)}$, and since $k$ is an integer we have $k\geq \lceil \sqrt{\chi(G)}\rceil$.

    Let us now show that $k \leq \lceil \sqrt{\chi(G)}\rceil$. For this, let $\ell$ be the smallest square number with $\ell \geq \chi(G)$, so we want to show that $k\leq \sqrt{\ell}$.
    Suppose that $\psi:V(G)\rightarrow [\ell]$ is a proper $\ell$-colouring of $G$.
    We take an arbitrary bijection $\rho$ between the set  $\{1,2,\ldots{},\ell\}$ and the set $\{(1,1),\ldots{},(1,\sqrt{\ell}),(2,1),\ldots{},(2,\sqrt{\ell}),\ldots{}(\sqrt{\ell},\sqrt{\ell})\}$, and assign the colours $i,j$ to the vertices $v_{in},v_{out}$ if $\rho(\psi(v))=(i,j)$. 
    To see that $\phi$ is an acyclic dicolouring of $D(G)$, it suffices to consider induced directed cycles of $D(G)$ which, by the remark above, are all 4-cycles of the form $C=(u_{in},u_{out},v_{in},v_{out})$ for some pair $u,v\in V$, and such a cycle corresponds to the edge $uv$ of $G$. As $\psi$ is a proper colouring of $G$, it follows from the definition of $\phi$ above that $\phi$ cannot assign the same colour to all four vertices of $C$, so each colour class is indeed acyclic. Suppose finally that $C$ alternates between $\phi^{-1}(i)$ and $\phi^{-1}(j)$. Then the definition of $\phi$ implies that  we have $\psi(u)=\psi(v)$, contradicting that $\psi$ is a proper colouring of $G$.
\end{proof}

It follows from \Cref{prop:chi_DG} that there exist bipartite tournaments ({\it i.e.} orientations of bipartite complete graphs) with arbitrarily high acyclic dichromatic number, and trivially all bipartite digraphs have dichromatic number at most 2. The same thing holds for tournaments by completing both parts of the partition into a transitive tournament. A {\bf transitive tournament} is an acyclic tournament, and the unique (up to isomorphism) transitive tournament of order $n$ is denoted $TT_n$.

\begin{proposition}
    \label{prop:dic2_dictk}
    For every integer $k$ there exists a tournament $T$ with $\dic{(T)=2}$,  $\dict{(T)}\geq k$ and such that $T$ has a spanning acyclic bipartite subdigraph.
\end{proposition}
\begin{proof}
    Let $B=(V_1,V_2,A)$ be a bipartite tournament with $\dict{(B)}\geq k$, which exists by \Cref{prop:chi_DG} (take for instance $D(K_{k^2})$), and form the tournament $T$ from two copies $B_1,B_2$ of $B$ with vertex sets $V_{1,1}\cup{}V_{1,2},V_{2,1}\cup{}V_{2,2}$ respectively by adding the arcs of a transitive tournament inside each $V_{i,j}$ and all arcs from $V_{1,1}\cup{}V_{1,2}$ to $V_{2,1}\cup{}V_{2,2}$.
    The partition $W_1=V_{1,1}\cup{}V_{2,1}$, $W_2=V_{2,1}\cup{}V_{2,2}$ shows that $\dic(T)=2$ and the fact that $B$ is a subdigraph of $T$ shows that $\dict{(T)}\geq k$.
    Moreover, $T$ has a spanning acyclic bipartite tournament induced by the arcs between $V_{1,1}\cup{}V_{1,2}$ and $V_{2,1}\cup{}V_{2,2}$.
\end{proof}

This shows that having dichromatic number 2 and a spanning acyclic bipartite subdigraph does not imply that the acyclic dichromatic number is bounded, even for tournaments. We pose the following problem.

\begin{problem}
    What is the complexity of deciding for given numbers $k_1,k_2\geq 1$ whether a tournament $T$ contains a spanning acyclic bipartite tournament $B=(V_1,V_2,A)$ where $|V_i|\geq k_i$ for $i=1,2$?
\end{problem}

A number of problems concerning the existence of spanning bipartite subdigraphs have been considered in \cite{bangJGT92}.

\section{Acyclic heroes of tournaments}\label{sec:acyclichero}

Let $T$ and $H$ be two tournaments. We say that $T$ is {\bf $H$-free} if $T$ does not contain any copy of $H$ as a subtournament. The class of $H$-free tournaments is denoted by $\forb{H}$. 
The tournament $H$ is a {\bf hero} if there exists a constant $c_{H}$ such that every tournament $T\in \forb{H}$ has dichromatic number at most $c_{H}$. 
Berger {\it et al.}~\cite{bergerJCT103} obtained an exact characterisation of heroes, for which we first need a few specific definitions. 

Let $H_1,H_2,H_3$ be three tournaments. We denote by $H_1\Ra H_2$ the tournament obtained from disjoint copies of $H_1$ and $H_2$ by putting all the arcs from $V(H_1)$ to $V(H_2)$. We further denote by $\Delta(H_1,H_2,H_3)$ the tournament obtained from disjoint copies of $H_1,H_2,H_3$ by putting all the arcs from  $V(H_1)$ to $V(H_2)$, from $V(H_2)$ to $V(H_3)$, and from $V(H_3)$ to $V(H_1)$.
For the sake of better readability, in the $\Delta$ operator, we write $k$ for the transitive tournament $TT_k$, so for instance $\Delta(H,1,k)$ denotes $\Delta(H,TT_1,TT_k)$.

We let $\mathcal{H}$ be the smallest class of tournaments such that:
\begin{enumerate}[label=$(\roman*)$]
    \item $\emptyset, TT_1 \in \mathcal{H}$,
    \item for every $H_1,H_2 \in \mathcal{H}$, the tournament $H_1 \Ra H_2$ belongs to $\mathcal{H}$, and
    \item for every $H\in \mathcal{H}$ and $k\geq 1$, both tournaments $\Delta(H,k,1)$ and $\Delta(H,1,k)$ belong to $\mathcal{H}$.
\end{enumerate}

It turns out that $\mathcal{H}$ is precisely the class of heroes.

\begin{theorem}[Berger {\it et al.}~\cite{bergerJCT103}]
    \label{thm:berger}
    A tournament $H$ is a hero if and only if $H\in \mathcal{H}$.
\end{theorem}

Similarly, a tournament $H$ is an {\bf acyclic hero} if there exists an absolute constant $c_{H}$ such that every tournament $T\in \forb{H}$ has acyclic dichromatic number at most $c_{H}$. 
Since $\dict(T) \leq \dic(T)$ holds for every tournament $T$, observe that, if $H$ is an acyclic hero, then it is also a hero, but the converse does not hold in general.

For every $k$, let $H_k$ be the tournament $TT_k \Ra (\Delta(k,1,1) \Ra TT_k)$. Let $\mathcal{AH}$ be the class of all tournaments appearing as a subtournament of $H_k$ for some integer $k$. We pose the following conjecture.

\begin{conjecture}
    \label{conj:heroes}
    A tournament $H$ is an acyclic hero if and only if $H\in \mathcal{AH}$.
\end{conjecture}

As evidences for the conjecture, we prove the following partial results in the next subsections.

\begin{theorem}
    \label{thm:non_acyclic_heroes}
    If $H$ is an acyclic hero, then $H\in \mathcal{AH}$.
\end{theorem}

\begin{theorem}
    \label{thm:heroes:C3_TT1}
    The tournament $\vec{C_3}\Ra TT_1$ is an acyclic hero.
\end{theorem}

\begin{theorem}
    \label{thm:heroes:Deltak11}
    For every integer $k$, $\Delta(k,1,1)$ is an acyclic hero. 
\end{theorem}

With Theorem~\ref{thm:heroes:C3_TT1} in hands, the following is a natural direction towards the resolution of Conjecture~\ref{conj:heroes}. Harutyunyan, Le, Thomassé, and Wu~\cite{harutyunyanJTCB138} proved the following ``local to global'' property for tournaments.

\begin{theorem}[Harutyunyan {\it et al.}~\cite{harutyunyanJTCB138}]
    There exists a function $f\colon \mathbb{N} \to \mathbb{N}$ such that every tournament $T$ satisfies
    \[
        \dic(T) \leq \max_{v\in V(T)}  f\Bigg( \dic\Big(T[N^+(v)]\Big) \Bigg).
    \]
\end{theorem}

We believe that the same property holds for the acyclic dichromatic number, and pose the following conjecture.

\begin{conjecture}
    \label{conj:local_to_global}
    There exists a function $f\colon \mathbb{N} \to \mathbb{N}$ such that every tournament $T$ satisfies
    \[
        \dict(T) \leq \max_{v\in V(T)}  f\Bigg( \dict\Big(T[N^+(v)]\Big) \Bigg).
    \]
\end{conjecture}

An easy inductive argument shows that, if true, Conjecture~\ref{conj:local_to_global} together with Theorem~\ref{thm:heroes:Deltak11} imply Conjecture~\ref{conj:heroes}.

\subsection{Proof of Theorem~\ref{thm:non_acyclic_heroes}}

\begin{lemma}
    \label{lemma:non_acyclic_heroes}
    The tournaments $\vec{C}_3 \Ra \vec{C}_3$, $\Delta(\vec{C_3},1,1)$, and $\Delta(2,2,1)$ are not acyclic heroes.
\end{lemma}
\begin{proof}
    Let us fix an arbitrary integer $k$, and exhibit a specific tournament $T$ with $\dict(T)\geq k$ that does not contain any of the aforementioned tournaments.

    Let $T$ be the tournament obtained as follows. We start from two transitive tournaments on $k^2$ vertices with acyclic orderings $u_1,\dots,u_{k^2}$ and $v_1,\dots,v_{k^2}$ respectively. Then, we add all arcs of the perfect matching $M=\{u_iv_i : i\in [1,k^2]\}$ and all arcs in $\{v_ju_i : i,j\in [1,k^2], i\neq j\}$. The following sequence of claims imply the lemma.

    \begin{claim}
        $\dict(T) \geq k$.
    \end{claim}
    \begin{proofclaim}
        This follows from the observation that $T$ contains $D(K_{k^2})$, which has total dichromatic number $k$ by \Cref{prop:chi_DG}.
    \end{proofclaim}
    \begin{claim}
        $T\in \forb{\vec{C_3} \Ra \vec{C_3}}$.
    \end{claim}
    \begin{proofclaim}
        Note that, in $T$, removing the arcs of $M$ yields an acyclic oriented graph. Therefore, every directed triangle contains at least one arc of $M$.
        Now assume for a contradiction that $T$ contains two disjoint directed triangles $C_1,C_2$ and all arcs from $V(C_1)$ to $V(C_2)$. By the remark above, there exists $1\leq i\neq j \leq k^2$ such that $\{u_i,v_i\} \subseteq V(C_1)$ and $\{u_j,v_j\} \subseteq V(C_2)$. By construction, $T$ contains the arc $v_ju_i$, a contradiction.
    \end{proofclaim}
    
    \begin{claim}
        $T\in \forb{\Delta(\vec{C_3},1,1)}$.
    \end{claim}
    \begin{proofclaim}
        Let us fix a directed triangle $C$. Again, $C$ has to use an arc of $M$, so let us denote, without loss of generality, $C=(u_i,v_i,x)$ for some $i\in [1,k^2]$ and some $x\in N^+(v_i) \cap N^-(u_i)$. Let $C^-$ be the set of vertices dominating $V(C)$, that is 
        \[
            C^- = N^-(u_i) \cap N^-(v_i) \cap N^-(x).
        \]
        Similarly, let $C^+$ be the set of vertices being dominated by $V(C)$. It follows from the definition of $T$ and the fact that $C^-\subseteq N^-(u_i) \cap N^-(v_i)$ that
        \[
        C^- \subseteq \{v_j : j\in [1,i-1]\},
        \]
        and that 
        \[
        C^+ \subseteq \{u_j : j\in [i+1,k^2]\},
        \]
        so indeed $C^-$ dominates $C^+$, implying that $T$ does not contain any copy of $\Delta(\vec{C_3},1,1)$.
    \end{proofclaim}
    \begin{claim}
        $T\in \forb{\Delta(2,2,1)}$.
    \end{claim}
    \begin{proofclaim}
        Observe first that, by construction, $T$ satisfies the following property:
          \begin{equation}
          \label{eqn:prop_star}
          \textit{For every $1\leq i<j\leq k^2$, both arcs $u_iu_j$ and $v_iv_j$ belong to exactly one directed triangle.}\tag{$\star$}
          \end{equation}
          Now assume for a contradiction that $T$ contains $\Delta(2,2,1)$, so there exist pairwise disjoint sets of vertices $\{x\}, Y=\{y_1,y_2\}, Z=\{z_1,z_2\}$ such that $x$ dominates $Y$, $Y$ dominates $Z$, and $Z$ dominates $x$. Observe that every arc in $\{xy_1,xy_2, z_1x,z_2x\}$ belongs to two distinct triangles. By~\eqref{eqn:prop_star}, it follows that each of these arcs belongs to the bipartite digraph $T[U,V]$ where $U=\{u_1,\dots,u_{k^2}\}$ and $V=\{v_1,\dots,v_{k^2}\}$.
          
          In particular, $x$ is a vertex with in-degree and out-degree at least $2$ in $T[U,V]$. This is a contradiction, as every vertex $u\in U$ has out-degree $1$ in $T[U,V]$, and every vertex $v\in V$ has in-degree $1$ in $T[U,V]$. The claim follows.
    \end{proofclaim}
    The lemma follows from the combination of the four previous claims.
\end{proof}

\begin{proof}[Proof of Theorem~\ref{thm:non_acyclic_heroes}]
    We prove the statement by induction on the order of $H$, the result being trivial when $|V(H)|\leq 1$. Suppose that $H$ is an acyclic hero of order at least $2$. In particular, $H$ is a hero, and by Theorem~\ref{thm:berger} one of the following holds:
    \begin{enumerate}[label=(\roman*)]
        \item there exist two smaller heroes $H_1, H_2$ such that $H$ is isomorphic to $H_1 \Ra H_2$, or
        \item there exist a smaller hero $H'$ and an integer $k$ such that $H$ is isomorphic to $\Delta(H',k,1)$ or $\Delta(H',1,k)$.
    \end{enumerate}
    
    In the first case, since $H_1$ and $H_2$ are subtournaments of $H$, they are both acyclic heroes, and by induction they both belong to $\mathcal{AH}$. If both of them contain $\vec{C_3}$, then $H$ contains $\vec{C_3} \Ra \vec{C_3}$, a contradiction to Lemma~\ref{lemma:non_acyclic_heroes}. Therefore, there exists $k$ such that one  of $H_1,H_2$ is isomorphic to $TT_k$ and the other one is isomorphic to a subtournament of $H_k=TT_k \Ra (\Delta(k,1,1) \Ra TT_k)$, hence implying that $H$ is isomorphic to a subtournament of $H_{2k}$.

    In the second case, if $H'$ contains $\vec{C_3}$, then $H$ contains $\Delta(\vec{C_3},1,1)$, a contradiction to Lemma~\ref{lemma:non_acyclic_heroes}. It follows that $H'$ is isomorphic to $TT_{k'}$ for $k'=|V(H')|$. If $\min(k,k')\geq 2$ then $H'$ contains $\Delta(2,2,1)$, a contradiction to Lemma~\ref{lemma:non_acyclic_heroes}. Therefore, $\min(k,k')=1$ and $H$ is isomorphic to $\Delta(\max(k,k'),1,1)$. The result follows.
\end{proof}

\subsection{Proof of Theorem~\ref{thm:heroes:C3_TT1}}

The proof of Theorem~\ref{thm:heroes:C3_TT1} relies on the particular structure of prime heroes. Given a tournament $T$, a {\bf homogeneous set} is a set of vertices $X$ such that, for every vertex $u\in V(T)\setminus X$, either $X\subseteq N^-(u)$ or $X\subseteq N^+(u)$. The set $X$ is further {\bf nontrivial} if $2\leq |X| \leq |V(T)|-1$.
A tournament is {\bf prime} if it does not contain any nontrivial homogeneous set. 

For an odd integer $n$, we define $R_n$ as the tournament with vertex set $V(R_n) = \{v_0,\dots,v_{n-1}\}$ and arc set
\[
A(R_n) = \left\{v_iv_{i+j} : 0\leq i \leq n-1 \mbox{ and } 1\leq j \leq \tfrac{n-1}{2} \right\},
\]
where indices are taken modulo $n$. Liu~\cite[Theorem~3.8]{liuArXiv2012} obtained the following result, from which we derive Theorem~\ref{thm:heroes:C3_TT1}.

\begin{theorem}[Liu~\cite{liuArXiv2012}]
    \label{thm:liu}
    Let $T$ be a prime tournament. If $n=|V(T)|\geq 3$ and $T\in \forb{\vec{C_3}\Ra TT_1}$, then $n$ is odd and $T$ is isomorphic to $R_n$.
\end{theorem}

\begin{proof}[Proof of Theorem~\ref{thm:heroes:C3_TT1}]
    We prove the stronger result that every tournament $T\in \forb{\vec{C_3}\Ra TT_1}$ satisfies $\dict(T)\leq 2$.

    For the sake of contradiction, let $T$ be a minimum counterexample, so $\dict(T) \ge 3$ and $\dict(T')\leq 2$ for every tournament $T'\in \forb{\vec{C_3}\Ra TT_1}$ of smaller order. By minimality, $T$ is strongly connected, otherwise we can $2$-colour each strongly connected component independently.

    We further claim that $T$ must be prime. Otherwise, let $X$ be a non-trivial homogeneous set. If $X$ contains a directed triangle, since $T$ is connected then $T$ contains a vertex $u$ with $X\subseteq N^-(u)$, and $T$ contains a subtournament isomorphic to $\vec{C_3}\Ra TT_1$, a contradiction. Hence, $T[X]$ is transitive. Let $T'$ be the tournament obtained from $T$ by contracting $X$ into a single vertex $x^\star$. By minimality of $T$, $T'$ admits an acyclic $2$-dicolouring $\phi$. It is straightforward to check that $\phi$ extends from $T'$ to $T$ by giving colour $\phi(x^\star)$ to all the vertices in $X$. 
    
    Hence, $T$ is prime, $|V(T)| \geq 3$, and $T\in \forb{\vec{C_3}\Ra TT_1}$. By Theorem~\ref{thm:liu}, $n=|V(T)|$ is odd and $T$ is isomorphic to $R_n$. It is straightforward to check that giving vertices $v_0,\dots,v_{(n-1)/2}$ colour $1$ and vertices $v_{(n+1)/2},\dots,v_{n-1}$ colour $2$ yields an acyclic $2$-dicolouring of $R_n$, hence yielding the contradiction.
\end{proof}

\subsection{Proof of Theorem~\ref{thm:heroes:Deltak11}}

We finally show that $\Delta(k,1,1)$ is an acyclic hero. We first show that it is an acyclic hero in the class of $2$-dicolourable tournaments, and then combine it with the fact that it is a hero to obtain the result.
We make use of the following celebrated result, obtained in 1964 by Erd\H{o}s and Moser~\cite{erdos1964}.
\begin{theorem}
    \label{thm:erdos_moser}
    For every integer $k$, every tournament of order at least $2^{k-1}$ contains $TT_k$.
\end{theorem}

\begin{lemma}
    \label{lemma:bip_Delta_k11}
    For every $k\geq 1$ and $T\in \forb{\Delta(k,1,1)}$, if $T$ satisfies $\dic(T)\leq 2$ then $\dict(T) \leq 2\cdot4^k$.
\end{lemma}
\begin{proof}
    We prove by induction that, for any tournament $T\in \forb{\Delta(k,1,1)}$ and any partition of $V(T)$ into two sets $(V_1,V_2)$ with $T[V_1]$ and $T[V_2]$ being acyclic, there exists an acyclic dicolouring of $T$ using at most $4^k$ colours on $V_1$ and at most $4^k$ colours on $V_2$, and that these two sets of colours are disjoint.

    Observe first that if one of $V_1$ or $V_2$ is empty, then the result is trivial, as we may colour every vertex of the other set with the same colour.
    Therefore, assume that $V_1\neq \emptyset$ and $V_2\neq \emptyset$. Let $t_1$ and $t_2$ be the sinks of $T[V_1]$ and $T[V_2]$, respectively. We assume by symmetry that the arc between $t_1$ and $t_2$ goes from $t_2$ to $t_1$. 

    Let $T' = T-t_1$, and assume by induction that $\phi$ is an acyclic dicolouring of $T'$, where $\phi(V_1\setminus \{t_1\}) \subseteq \{1,\dots,4^k\}$ and $\phi(V_2) \subseteq \{4^k+1,\dots,2\cdot4^k\}$. We claim that $\phi$ can be extended to $T$ by choosing for $t_1$ a colour from $\{1,\dots, 4^k\}$ that is not given by $\phi$ to any vertex in 
    \[
        \bigcup_{u\in N^+(t_1)} N^+(u).
    \]

    Assume first that such a colour exists, and let us briefly justify that the obtained extension of $\phi$ is indeed an acyclic dicolouring. If this is not the case, then there exists a directed cycle $\Cscr$ which is either monochromatic or alternating between two colour classes. Among all such cycles, we let $\Cscr$ be a shortest one, so it has length either $3$ or $4$ (as $T$ is a tournament). By the induction hypothesis, $\Cscr$ contains $t_1$. Moreover, since $T[V_1]$ is acyclic, $\Cscr$ contains at least one vertex of $V_2$, and because vertices in $V_1$ and $V_2$ use distinct sets of colours, $\Cscr$ is necessarily a directed $4$-cycle $(t_1,x,y,z)$ alternating between two colour classes. We conclude that $\phi(t_1)=\phi(y)$, a contradiction, since $y\in \bigcup_{u\in N^+(t_1)} N^+(u)$.

    We now argue that such a colour exists. Recall that $N^+(t_1) \subseteq V_2$, and that $N^+(t_1) \subseteq N^-(t_2)$, since $t_1$ and $t_2$ are sinks in $T[V_1]$ and $T[V_2]$ respectively. Therefore, we have $d^+(t_1) \leq 2^{k-1}$, for otherwise, by Theorem~\ref{thm:erdos_moser}, $T[N^+(t_1)\cap N^-(t_2)]$ contains a copy of $TT_k$, contradicting the assumption that $T\in \forb{\Delta(k,1,1)}$.

    Now, let us fix a vertex $u\in N^+(t_1)$. The same argument, relying on the fact that $t_1$ is a sink in $T[V_1]$, shows that $|N^+(u) \cap V_1| \leq 2^{k-1}$. It follows that $\bigcup_{u\in N^+(t_1)} N^+(u)$ contains at most $2^{2k-2} \le 4^k-1$ vertices of $V_1$. Since vertices of $V_2$ use colours distinct from $\{1,\dots,4^k\}$, it follows that one colour is available for $t_1$, showing the result.
\end{proof}

\begin{lemma}
    \label{lemma:heroes_dic2}
    Let $H$ be a hero such that, for some $c \in \mathbb{N}$, every $H$-free tournament $T$ with $\dic(T)\leq 2$ satisfies $\dict(T) \leq c$. Then $H$ is an acyclic hero.
\end{lemma}
\begin{proof}
    Let $T$ be an arbitrary $H$-free tournament, and let us show that it has bounded acyclic dichromatic number. 

    As $H$ is a hero, there exists $k$ such that $\dic(T)\leq k$. We fix $V_1,\dots,V_{k}$ a partition of $V(T)$ into transitive tournaments. For every $1\leq i<j\leq k$, we let $\phi_{i,j}$ be an acyclic dicolouring of $T[V_i \cup V_j]$ using at most $c$ colours, which exists by assumption.
    For every $i\in [1,k]$ and every vertex $u\in V_i$, we let 
    \[
        \psi(u) = \Big(\phi_{1,i}(u), \phi_{2,i}(u), \dots, \phi_{i-1,i}(u), 0, \phi_{i,i+1}(u), \dots,\phi_{i,k}(u)\Big).
    \]
    For any vertex $u$ and integer $i\in [1,k]$, we further denote by $\psi_i(u)$ the $i^{\rm th}$ coordinate of $\psi(u)$.
    Note that $\psi$ is a colouring of $T$ using at most $(c+1)^{k}$ colours (since every coordinate of $\psi(u)$ is an integer of $[0,c]$). It remains to show that it is indeed an acyclic dicolouring. 

    Assume for a contradiction that $T$, coloured with $\psi$, contains a directed cycle $\Cscr$ which is either monochromatic or alternating between two colour classes of $\psi$. Among all such cycles, we choose $\Cscr$ of minimum length.
    
    Note that, by definition, $\psi_i(u) = 0$ if and only if $u\in V_i$. In particular, each colour class of $\psi$ is included in some $V_i$. It follows that $\Cscr$ cannot be monochromatic, so it is alternating between $V_i$ and $V_j$ for some indices $1\leq i < j\leq k$, and $\Cscr$ has length $4$. Let us denote $\Cscr=(u,x,v,y)$, where $u,v\in V_i$ and $x,y\in V_j$.

    Since $\Cscr$ alternates between two colour classes of $\psi$, in particular we have $\psi_j(u)=\psi_j(v)$, which by definition of $\psi$ implies $\phi_{i,j}(u) = \phi_{i,j}(v)$. Similarly, we have $\phi_{i,j}(x) = \phi_{i,j}(y)$. Therefore, $\Cscr$ is a directed cycle in $T[V_i\cup V_j]$ alternating between two colour classes of $\phi_{i,j}$, a contradiction. The result follows.
\end{proof}

Theorem~\ref{thm:heroes:Deltak11} follows directly from the combination of Lemmas~\ref{lemma:bip_Delta_k11} and~\ref{lemma:heroes_dic2}.

\section{Complexity} 
\label{sec:complexity}

In this section, we study the complexity of the following problem, defined for every fixed integer~$k$.

\defproblem{
    \sc Acyclic $k$-Dicolourability
}
{
    An oriented graph $D$.
}{
    Does $D$ admit an acyclic $k$-dicolouring?
}

Recall that, by \Cref{prop:chi_DG}, for any graph $G$, one can construct in polynomial time the oriented graph $D(G)$ satisfying
$\dict(D(G)) = \lceil \sqrt{\chi(G)}\rceil$. Recall furthermore that $D(G)$ is a bipartite digraph.
Since, for every integer $k\geq 2$, it is NP-complete to decide whether a graph is $k^2$-colourable, we directly obtain the following.

\begin{corollary}
\label{cor:NPCbipartite}
    For every integer $k\geq 2$, {\sc Acyclic $k$-Dicolourability} is NP-complete even when restricted to bipartite digraphs.
\end{corollary}

In Section~\ref{sec:split}, we give another reduction showing that {\sc Acyclic $k$-Dicolourability} is NP-hard, for every integer $k\geq 3$, even on the class of split digraphs. In Section~\ref{subsec:complexity:tournaments}, we provide a polynomial-time algorithm for {\sc Acyclic $2$-Dicolourability} on tournaments.

\subsection{Complexity for split digraphs}
\label{sec:split}

A digraph $D=(V,A)$ is a {\bf split digraph} if we can partition its vertices into two sets $V_1,V_2$ so that $D[V_1]$ has no arcs and $D[V_2]$ is a semicomplete digraph. We call $V_1,V_2$ a {\bf split partition of $V(D)$.} A split digraph $D=(V,A)$ is {\bf semicomplete} if every vertex in $V_1$ is adjacent to every vertex in $V_2$.

\begin{theorem}
    \label{thm:NP_hardness_split}
    For every fixed $k\geq 3$, {\sc Acyclic $k$-Dicolouring} is NP-complete even when restricted to split digraphs.
\end{theorem}

However, the complexity of the problem remains open when $k=2$ or when the input is a complete split digraph.

\begin{problem}
What is the complexity of deciding whether a split digraph $D$ has $\dict(D)=2$?
\end{problem}

\begin{problem}
What is the complexity of deciding whether a semicomplete split digraph $D$ has $\dict(D)=2$?
\end{problem}

Note that a semicomplete split digraph is a special case of a semicomplete multipartite digraph in which all but one partite set has size one. Hence, proving the above polynomial would be a step towards determining the complexity for multipartite tournaments.

The remainder of this section is devoted to the proof of Theorem~\ref{thm:NP_hardness_split}. We first prove the hardness for $k=3$.

\begin{lemma}
    \label{lemma:NP_hardness_split_3}
    {\sc Acyclic $3$-Dicolouring} is NP-hard even when restricted to split digraphs.
\end{lemma}
\begin{proof}
    We reduce from {\sc $2$-Dicolourability}, which is well-known to be NP-hard (see~\cite{bokalJGT46}). Let $D$ be any instance of this problem.
    We assume that $D$ has minimum in-degree at least $1$, for otherwise we repeatedly remove sources from $D$ to obtain an equivalent instance.
    
    We build $D^\star$ as follows. 
    To every vertex $v\in V(D)$ we associate a vertex $x_v$ in $D^\star$, and to every arc $a\in A(D)$ we associate a vertex $x_a$. Furthermore, for every arc $a=uv\in A(D)$, we add the arcs $x_ux_a,x_ax_v$. We fix an arbitrary ordering $a_1,\dots ,a_{m}$ of $A(D)$, and add to $D^\star$ all the arcs $x_{a_i}x_{a_j}$ with $j>i$, so $X_A = \{x_a:a\in A(D)\}$ induces a transitive tournament. We further add a new vertex $x_{b}$ dominating $X_A$ and being dominated by $X_V=\{x_v:v\in V(D)\}$
    We finally add the six vertices in 
    $Y=\{y_1,y_2,y_3\}$ and $Z=\{z_1,z_2,z_3\}$, together with the arcs in
    \[
        \{y_1y_2, y_2y_3, y_3y_1\} \cup \{z_1z_2, z_2z_3, z_3z_1\} \cup \{ x_ay, yz,zx_a : a\in A(D)\cup\{b\}, y\in Y, z\in Z\},
    \]
    see Figure~\ref{fig:split_reduction} for an illustration. 

    \begin{figure}[ht]
        \centering
        \begin{tikzpicture}
        \def\L{1.2}
        \def\H{1}
        \node[svertex, orange, label=left:$t$] (t) at (0,0) {};
        \node[svertex, g-blue, label=left:$u$] (u) at (0,-\H) {};
        \node[svertex, g-blue, label=left:$v$] (v) at (0,-2*\H) {};
        \node[svertex, orange, label=left:$w$] (w) at (0,-3*\H) {};

        \draw[sarc] (t) to (u);
        \draw[sarc] (u) to (v);
        \draw[sarc] (v) to (w);
        
        \draw[sarc] (v) to[out=60, in=-60] (t);
        \draw[sarc, ultra thick,white] (w) to[out=60, in=-60] (u);
        \draw[sarc] (w) to[out=60, in=-60] (u);

        \begin{scope}[xshift=4cm]
            \node[svertex, orange, label=left:$x_t$] (xt) at (0,0) {};
            \node[svertex, g-blue, label=left:$x_u$] (xu) at (0,-\H) {};
            \node[svertex, g-blue, label=left:$x_v$] (xv) at (0,-2*\H) {};
            \node[svertex, orange, label=left:$x_w$] (xw) at (0,-3*\H) {};
            
            \node[svertex, g-green] (xtu) at (\L,\H) {};
            \node[svertex, g-green] (xvt) at (\L,0) {};
            \node[svertex, g-green] (xuv) at (\L,-\H) {};
            \node[svertex, g-green] (xwu) at (\L,-2*\H) {};
            \node[svertex, g-green] (xvw) at (\L,-3*\H) {};
            \node[svertex, g-green, label=right:$x_b$] (xb) at (\L,-4*\H) {};
            
            \node[] (blank) at (1.8*\L,0) {};
            
            \node[rectangle, rounded corners, draw=gray!70, inner xsep=2pt, inner ysep=5pt, fill=none, fit=(xt) (xw), label=above:{\color{gray} \footnotesize $X_V$}] (XV) {};
            
            \node[rectangle, rounded corners, draw=gray!70, inner sep=5pt, fill=none, fit=(xtu) (xb) (blank), label=above:{\color{gray} \footnotesize $X_A \cup \{x_b\}$}] (XA) {};

            \draw[sarc] (xt) to (xtu);
            \draw[sarc] (xtu) to (xu);
            
            \draw[sarc] (xv) to (xvt);
            \draw[sarc, ultra thick, white] (xvt) to (xt);
            \draw[sarc] (xvt) to (xt);
            
            \draw[sarc, ultra thick, white] (xu) to (xuv);
            \draw[sarc] (xu) to (xuv);
            \draw[sarc] (xuv) to (xv);
            
            \draw[sarc] (xw) to (xwu);
            \draw[sarc, ultra thick, white] (xwu) to (xu);
            \draw[sarc] (xwu) to (xu);
            
            \draw[sarc, ultra thick, white] (xv) to (xvw);
            \draw[sarc] (xv) to (xvw);
            \draw[sarc] (xvw) to (xw);

            \draw[sarc] (XV.-90) to[out=-90, in=180] (xb);
            \draw[sarc, gray!80] (xb) to (xvw);

            \draw[sarc, gray!80] (xvt) to (xtu);
            \draw[sarc, gray!80] (xuv) to (xvt);
            \draw[sarc, gray!80] (xwu) to (xuv);
            \draw[sarc, gray!80] (xvw) to (xwu);
            
            \draw[sarc, gray!80] (xuv) to[out=45, in=-45]  (xtu);
            \draw[sarc, gray!80] (xwu) to[out=45, in=-45]  (xvt);
            \draw[sarc, gray!80] (xwu) to[out=45, in=-45]  (xtu);
            \draw[sarc, gray!80] (xvw) to[out=45, in=-45] (xuv);
            \draw[sarc, gray!80] (xvw) to[out=45, in=-45]  (xvt);
            \draw[sarc, gray!80] (xvw) to[out=45, in=-45]  (xtu);

            \draw[sarc, gray!80] (xb) to[out=45, in=-45]  (xwu);
            \draw[sarc, gray!80] (xb) to[out=45, in=-45] (xuv);
            \draw[sarc, gray!80] (xb) to[out=45, in=-45]  (xvt);
            \draw[sarc, gray!80] (xb) to[out=45, in=-45]  (xtu);
            \begin{scope}[xshift=4.5cm, yshift=-0.2cm]
                \node[svertex, orange] (y1) at (90:0.5) {};
                \node[svertex, orange] (y2) at (-30:0.5) {};
                \node[svertex, g-blue] (y3) at (-150:0.5) {};
                \draw[sarc] (y1) to (y2);
                \draw[sarc] (y2) to (y3);
                \draw[sarc] (y3) to (y1);   
                \node[circle, rounded corners, draw=gray!70,  fill=none, fit=(y1) (y2) (y3), inner sep=0, label=right:{\color{gray} \footnotesize $Y$}] at (0,0) (Y) {};
            \end{scope}
            \begin{scope}[xshift=4.5cm, yshift=-2.8cm]
                \node[svertex, orange] (z1) at (90:0.5) {};
                \node[svertex, orange] (z2) at (-30:0.5) {};
                \node[svertex, g-blue] (z3) at (-150:0.5) {};
                \draw[sarc] (z1) to (z2);
                \draw[sarc] (z2) to (z3);
                \draw[sarc] (z3) to (z1);  
                \node[circle, rounded corners, draw=gray!70,  fill=none, fit=(z1) (z2) (z3), inner sep=0, label=right:{\color{gray} \footnotesize $Z$}] at (0,0) (Z) {};
            \end{scope}
            \draw[sarc] (Y) to (Z);
            \draw[sarc] (Z) to (XA);
            \draw[sarc] (XA) to (Y);
        \end{scope}
        \end{tikzpicture}
        \caption{A digraph $D$ together with one of its $2$-dicolourings $\phi$ (left), and the corresponding digraph $D^\star$ with one of its acyclic $3$-dicolourings (right).}
        \label{fig:split_reduction}
    \end{figure}
    
    By construction, $X_V$ is an independent set, and $D^\star - X_V$ is a tournament, so $D^\star$ is a split digraph. We now show that $\dic(D)\leq 2$ if and only if $\dict(D^\star)\leq 3$, implying the lemma.

    \begin{claim}
        If $\dic(D)\leq 2$ then $\dict(D^\star)\leq 3$. 
    \end{claim}
    \begin{proofclaim}
        Assume that $\phi\colon V(D)\to \{1,2\}$ is a $2$-dicolouring of $D$. Let $\phi^\star\colon V(D^\star) \to \{1,2,3\}$ be the colouring of $D^\star$ defined as follows:
        \[
            \phi^\star(x) = \left\{
            \begin{array}{ll}
                \phi(v) & \mbox{if } x = x_v \mbox{ with }v\in V(D) \\
                3 & \mbox{if } x \in X_A \cup \{x_b\} \\
                1 & \mbox{if } x \in \{y_1,z_1\} \\
                2 & \mbox{otherwise, that is if } x\in \{y_2,y_3,z_2,z_3\}.
            \end{array}
        \right.   
        \]
        We claim that $\phi^\star$ is an acyclic dicolouring of $D^\star$. Clearly, each colour class induces an acyclic digraph, and the digraph $D[V_1,V_2]$ is acyclic as well (where $V_i$ is the set of vertices coloured $i$ in $\phi^\star$). Therefore, if $\phi^\star$ is not an acyclic dicolouring, there exists a directed cycle $\Cscr$ alternating between $V_i$ and $V_3$ for some $i\in \{1,2\}$. Note that $\Cscr$ does not contain any vertex $z\in Z$, as every in-neighbour of $z$ is coloured with $\{1,2\}$. Similarly, $\Cscr$ does not contain any vertex in $Y$, and it does not contain $x_b$. Therefore, $\Cscr$ is included in the digraph $D[X_V,X_A]$, meaning that $\Cscr$ is of the form 
        \[
            \Cscr = (x_{u_1},x_{u_1u_2},x_{u_2},\ldots, x_{u_{\ell-1}u_\ell},x_{u_\ell},x_{u_\ell u_1}).
        \]
        But then $\Cscr'= (u_1,u_2,\dots,u_\ell)$ is a monochromatic directed cycle in $D$ coloured with $\phi$, a contradiction.
    \end{proofclaim}
    
    \begin{claim}
        If $\dict(D^\star)\leq 3$ then $\dic(D)\leq 2$. 
    \end{claim}
    \begin{proofclaim}
        Let $\phi^\star$ be an acyclic $3$-dicolouring of $D$ with colour classes $V_1,V_2,V_3$, and assume without loss of generality that $x_b\in V_3$.
        
        We first show that $X_A\cup \{x_b\} \subseteq V_3$. Assume that this is not the case, and assume by symmetry that $x_a\in V_1$ for some $a\in A(D)$. Since both $Y$ and $Z$ induce directed triangles, note that some colour appears on both $Y$ and $Z$. If colour $3$ belongs to $\phi(Y) \cap \phi(Z)$, then $(x_b,y,z)$ is a monochromatic directed triangle for some $y\in Y$ and $z\in Z$. Similarly, because of $x_a$ being coloured $1$, colour $1$ does not belong to $\phi(Y) \cap \phi(Z)$. Therefore, $\phi(Y) \cap \phi(Z) = \{2\}$.
        If there exists $y\in Y$ with $\phi(y) = 3$, then there exists $z\in Z$ with $\phi(z)=1$, and $(x_a,y,z,x_b)$ alternates between $V_1$ and $V_3$, a contradiction. Therefore, $\phi(Y) = \{1,2\}$ and $\phi(Z) = \{2,3\}$. Assume without loss of generality that $\phi(y_1) = 2$, $\phi(z_1)=3$ and $\phi(z_2)=2$, then $(x_b,y_1,z_1,z_2)$ alternates between $V_2$ and $V_3$, a contradiction. This shows $X_A\cup \{x_b\} \subseteq V_3$.  

        It follows that $V_3 \cap X_V = \emptyset$. Indeed, if some vertex $x_v\in X_V$ is coloured $3$, then for any entering arc $a$ of $v$ in $D$ (recall that $D$ has minimum in-degree at least $1$), then the 3-cycle $(x_b,x_a,x_v)$ is monochromatic, a contradiction. This allows us to define $\phi\colon V(D)\to \{1,2\}$, the $2$-colouring of $D$ where, for every $v\in V(D)$,
        \[
            \phi(v) = \phi^\star(x_v).
        \]
        We claim that $\phi$ is a $2$-dicolouring of $D$. Indeed, if $\Cscr= (u_1,u_2,\dots,u_\ell)$ is a monochromatic directed cycle in $D$ coloured with $\phi$, then
        \[
            \Cscr^\star= (x_{u_1},x_{u_1u_2},x_{u_2},\ldots, x_{u_{\ell-1}u_\ell},x_{u_\ell},x_{u_\ell u_1})
        \]
        is a directed cycle in $D^\star$ alternating between two colour classes of $\phi^\star$, a contradiction.
    \end{proofclaim}
    The lemma follows.
\end{proof}

We now derive \Cref{thm:NP_hardness_split}.

\begin{proof}[Proof of \Cref{thm:NP_hardness_split}]
    It is clear that {\sc Acyclic $k$-Dicolourability} is in NP. We prove by induction on $k\geq 3$ that it is NP-hard for every $k\geq 3$, even when restricted to split digraphs, the case $k=3$ corresponding exactly to \Cref{lemma:NP_hardness_split_3}.

    Let $D_k$ be any digraph and $T_k$ be a tournament with acyclic dichromatic number $k$, whose size depends only on $k$ ({\it e.g.} the tournament constructed in the proof of \Cref{prop:dic2_dictk}, which has $\dict \geq k$ so contains a subtournament, $T_k$, with $\dict(T_k)=k$).
    We let $D_{k+1}$ be the digraph obtained by taking disjoint copies $D_k$ and $T_k$ with one extra vertex $v$ and putting all arcs in
    \[
        \{vu,ut,tv : u\in V(D_k), t\in V(T_k)\}.
    \]
    Note that the size of $D_{k+1}$ is at most a polynomial in the size of $D_k$ ($k$ being fixed), and that if $D_k$ is a split digraph then so is $D_{k+1}$.
    The hardness of {\sc Acyclic $(k+1)$-Dicolourability} for split digraphs thus follows from the hardness of {\sc Acyclic $k$-Dicolourability} for split digraphs together with the following claim.  Note that we may assume that $D_k$ is strongly connected, as when computing the acyclic dicolourability of a digraph, it suffices to do so for each of the strong components (and take the maximum over all strong components).
    \begin{claim}
        We have $\dict(D_k) \leq k$ if and only if $\dict(D_{k+1}) \leq k+1$.
    \end{claim}
    \begin{proofclaim}
        Assume first that $\dict(D_{k}) \leq k$, and let $\psi$ be the $(k+1)$-colouring of $D_{k+1}$ obtained by taking an acyclic $k$-dicolouring of $D_k$, an acyclic $k$-dicolouring of $T_k$, and using colour $k+1$ for $v$. Assume for a contradiction that $\psi$ is not an acyclic dicolouring of $D_{k+1}$, that is $D_{k+1}$ contains a directed cycle $\Cscr$ which is either monochromatic or alternating between two colours. Then $\Cscr$ does not contain $v$, as $v$ is the only vertex coloured $k+1$. Hence, $\Cscr$ is included in a strongly connected component of $D_{k+1}-v$, so in particular $\Cscr$ is included either in $D_k$ or in $T_k$, a contradiction.

        Assume now that $\dict(D_{k+1}) \leq k+1$ and let $\psi$ be an acyclic $(k+1)$-dicolouring of $D_{k+1}$. We assume without loss of generality that $v$ is coloured $k+1$. If both $T_k$ and $D_k$ use colour $k+1$, then $D_{k+1}$ contains a monochromatic directed triangle, a contradiction. Therefore, one of them does not use colour $k+1$. If $D_k$ does not use colour $k+1$, then $\psi$ uses at most $k$ colours on $D_k$ and the result follows, so we assume that $(k+1) \notin \psi(V(T_k))$. 
Since $\dict(T_k)=k$, for every $i\in [1,k]$ there exists a vertex $v_i$ coloured $i$ in $T_k$. 
Let $\mathcal{U}$ be the set of vertices coloured $k+1$ in $D_k$. 
Clearly, we may assume that $\mathcal{U}$ has at least one entering arc in $D_k$, for otherwise $D_k$ would not be strongly connected.
Hence, let $xy\in A(D_k)$ be such that $\psi(x) \neq k+1$ and $\psi(y) = k+1$. Now there exists $z\in V(T_k)$ with $\psi(z) = \psi(x)$, and $(v,x,y,z)$ is a directed cycle alternating between colour classes $k+1$ and $\psi(x)$, a contradiction. 
    \end{proofclaim}
    The result follows.
\end{proof}

\subsection{A polynomial-time algorithm for {\sc Acyclic \texorpdfstring{$2$}{2}-Dicolouring} on tournaments}
\label{subsec:complexity:tournaments}

We proved in the last sections that deciding $\dict(D)\leq 2$ is NP-hard for general oriented graphs $D$. We now show that the situation is different for tournaments. This also draws a distinction with the {\sc $2$-Dicolourability} problem, which is well-known to be NP-hard even when restricted to tournaments~\cite{chenSIC37} (see also~\cite[Theorem~2.8.21]{bang2018} for a shorter proof).

Our proof is inspired by a proof due to Klingelhoefer and Newman~\cite{klingelhoeferSIDMA38} that a $10$-dicolouring of a tournament $T$ can be computed in polynomial time when $T$ is guaranteed to have dichromatic number $2$. Precisely, we make use of a decomposition of $2$-dicolourable light tournaments that they obtained. A tournament $T$ is {\bf light} if, for every arc $uv\in A(T)$, the subtournament $T[N^+(v)\cap N^{-}(u)]$ is acyclic. 

\begin{theorem}
    \label{thm:polytime_tour_2colours}
    There exists a polynomial-time algorithm for deciding whether a given tournament $T$ admits an acyclic $2$-dicolouring, and computing one when it exists.
\end{theorem}

The cases of more than two colours remain open. We believe that {\sc Acyclic $k$-Dicolouring}, when restricted to tournaments, is actually solvable in polynomial time for every fixed $k$. We pose the following stronger conjecture. A tournament $T$ is {\bf $k$-critical} if $\dict(T) = k$ and every proper subtournament $H$ of $T$ satisfies $\dict(H) \leq k-1$. 

\begin{conjecture}
    For every $k\in \mathbb{N}$, there is a finite number of $k$-critical tournaments.
\end{conjecture}

If the conjecture above is true, then in particular there is a trivial polynomial time algorithm for deciding $\dict(T)\leq k$, consisting of checking whether $T$ contains one of the $(k+1)$-critical tournaments. The conjecture remains open for every $k\geq 2$.

The remainder of this subsection is devoted to the proof of Theorem~\ref{thm:polytime_tour_2colours}.

\begin{proof}[Proof of Theorem~\ref{thm:polytime_tour_2colours}]
    Let $T$ be any tournament on $n\geq 3$ vertices (the algorithm being trivial when $n\leq 2$). Let us describe the claimed algorithm. Clearly, we may assume that $T$ is strongly connected, for otherwise we apply the following algorithm to each strongly connected component of $T$. We first justify that we can further assume that $T$ is light.

    \begin{claimlabeled}[thm:polytime_tour_2colours]
        If $\dict(T)\leq 2$ then $T$ is light.
        \label{claim:T_light}
    \end{claimlabeled}
    \begin{proofclaim}
        Assume that $\dict(T)\leq 2$ and let $\phi$ be an acyclic $2$-dicolouring of $T$. Assume for a contradiction that $T$ is not light, so there exists an arc $uv \in A(T)$ such that $T[N^+(v) \cap N^-(u)]$ contains a directed triangle $(x,y,z)$.

        Assume first that $\phi(u) = \phi(v)$. Since $(x,y,z)$ is a directed triangle, one of these three vertices, say $x$, uses colour $\phi(u)$. Therefore, $(v,x,u)$ is a monochromatic directed triangle, a contradiction.

        Assume now, without loss of generality, that $\phi(u) = 1$ and $\phi(v)=2$. Since $\phi$ uses both colours $1$ and $2$ on $\{x,y,z\}$ (otherwise $(x,y,z)$ is monochromatic), there exists an arc, say $xy$, such that $\phi(x) = 1$ and $\phi(y)=2$. Therefore, $(v,x,y,u)$ is a directed cycle alternating between colours $1$ and $2$, a contradiction to $\phi$ being an acyclic dicolouring.
    \end{proofclaim}
    
    Clearly, one can check in polynomial time whether $T$ is light. If it is not, then by \Cref{claim:T_light} we have $\dict(T) \geq 3$ and we can stop here.
    From now on, we thus assume that $T$ is light.

    An ordered pair $(s,t)$ of distinct vertices is {\bf nice} if there exists an acyclic $2$-dicolouring $\phi$ of $T$ such that $\phi(s) = \phi(t) = 1$ and, for every $u\in N^+(s) \cup N^-(t)$, $\phi(u) = 2$. Such a colouring $\phi$ is called an {\bf $(s,t)$-colouring} of $T$.
    \begin{claimlabeled}[thm:polytime_tour_2colours]
        \label{claim:existence_nice_pair}
        $\dict(T)\leq 2$ if and only if $T$ admits a nice pair $(s,t)$.
    \end{claimlabeled}
    \begin{proofclaim}
        By definition, if $T$ admits a nice pair then it admits an acyclic $2$-dicolouring, and in particular $\dict(T)\leq 2$.
        
        Conversely, assume that $\dict(T)\leq 2$ and let $\phi$ be an acyclic $2$-dicolouring of $T$. Let $V_1$  and $V_2$ be the vertices coloured respectively $1$ and $2$ in $\phi$, and assume without loss of generality that $|V_1|\geq 2$.
        Since $T[V_1]$ is transitive, let $s$ be its sink and $t$ be its source. By definition, $\phi$ is an acyclic $2$-dicolouring of $T$ such that $\phi(s)= \phi(t)=1$. Moreover, $N^+(s) \cup N^-(t) \subseteq V_2$, hence showing that $\phi$ is an $(s,t)$-colouring.
    \end{proofclaim}

    From now on, let us fix two arbitrary distinct vertices $s,t$. 
    In the remaining of the proof, we show that we can decide in polynomial time whether $(s,t)$ is a nice pair, and computing an $(s,t)$-colouring when it is.
    Since there exist at most $n(n-1)$ candidates $(s,t)$ for being a nice pair, the result then follows from \Cref{claim:existence_nice_pair}.

    We first compute a shortest directed $(s,t)$-path $P=(v_0,v_1,\dots,v_{\ell-1}, v_\ell)$, where $s=v_0$ and $t=v_\ell$ (recall that $T$ is strongly connected). We sequentially build a sequence $(X_0,\dots,X_\ell)$ of disjoint subsets of $V(T)$ where:
    \begin{align*}
        X_0 &= N^+(v_0) \setminus V(P)\\
        \text{and~~} X_i &= N^+(v_i) \setminus \Big(V(P) \cup \bigcup_{j=0}^{i-1}X_j\Big) \text{~~for every $i\in [1,\ell]$.}
    \end{align*}
    Furthermore, let $X_{\ell+1} = N^-(v_\ell)\setminus \Big(V(P) \cup \bigcup_{j=0}^{\ell}X_j\Big)$.
    The point is that now some structure is forced,
    see Figure~\ref{fig:nice_decomposition} for an illustration.

    \begin{figure}[ht]
        \begin{center}  
              \begin{tikzpicture}
                \def\L{1.8};
                \def\H{1.2};
                \node[svertex, label=left:{$s=v_0$}] (v0) at (0,0) {};
                \node[svertex, label=below:$v_1$] (v1) at (\L,0) {};
                \node[svertex, label=below:$v_2$] (v2) at (2*\L,0) {};
                \node[svertex, label=below:$v_3$] (v3) at (3*\L,0) {};
                \node[svertex, label=below:$v_4$] (v4) at (4*\L,0) {};
                \node[svertex, label=below:$v_{\ell-2}$] (v5) at (5*\L,0) {};
                \node[svertex, label=below:$v_{\ell-1}$] (v6) at (6*\L,0) {};
                \node[svertex, label=right:{$v_{\ell}=t$}] (v7) at (7*\L,0) {};

                \node[labelledvertex] (X0) at (0,-\H-0.5) {$X_0$};
                \node[labelledvertex] (X1) at (0.5*\L,\H) {$X_1$};
                \node[labelledvertex] (X2) at (1.5*\L,\H) {$X_2$};
                \node[labelledvertex] (X3) at (2.5*\L,\H) {$X_3$};
                \node[labelledvertex] (X4) at (3.5*\L,\H) {$X_4$};
                \node[labelledvertex] (X6) at (5.5*\L,\H) {$X_{\ell-1}$};
                \node[labelledvertex] (X7) at (6.5*\L,\H) {$X_\ell$};
                \node[labelledvertex] (X8) at (7*\L,-\H-0.5) {$X_{\ell+1}$};

                \foreach \i in {0,1,2,3,4,6,7}{
                    \draw[sarc] (v\i) to (X\i);
                }
                \draw[sarc] (X8) to (v7);
                \foreach \i in {0,1,2,3,5,6}{
                    \pgfmathtruncatemacro{\j}{\i+1};
                    \draw[sarc] (v\i) to (v\j);
                }
                \draw[sarc, thick, dotted] (v4) to (v5);
                \foreach \i in {1,2,3,4,6,7}{
                    \pgfmathtruncatemacro{\j}{\i-1};
                    \draw[sarc] (X\i) to (v\j);
                }
              \end{tikzpicture}
          \caption{The partition $(V(P),X_0,\dots,X_{\ell+1})$.}
          \label{fig:nice_decomposition}
        \end{center}
    \end{figure}
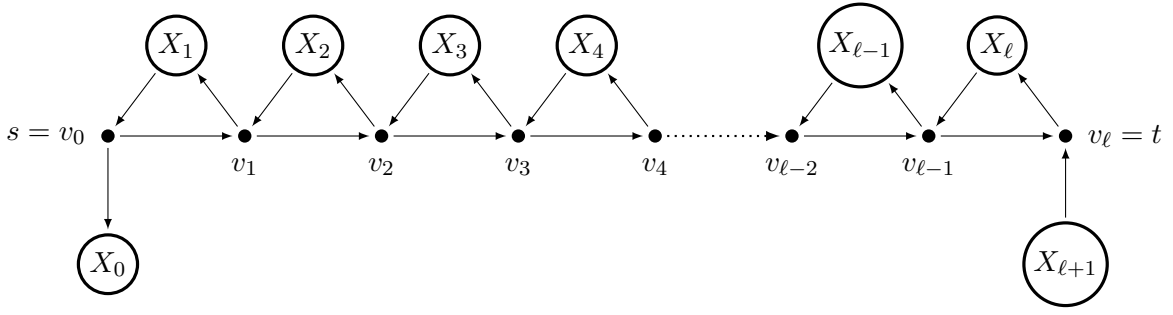
    
    \begin{claimlabeled}[thm:polytime_tour_2colours]
        \label{claim:structure_nice_decomposition}
        The following properties hold:
        \begin{enumerate}[label={\rm (\roman*)}]
            \item $X_0 \subseteq N^+(s)$ and $X_{\ell+1} \subseteq N^-(t)$;
            \label{claim:structure_nice_decomposition:enum:1}
            \item for every $i \in [1,\ell]$, $X_i \subseteq N^+(v_i) \cap N^-(v_{i-1})$;
            \label{claim:structure_nice_decomposition:enum:2}
            \item for every $i \in [1,\ell]$, $T[X_i]$ is acyclic; and
            \label{claim:structure_nice_decomposition:enum:3}
            \item $(V(P), X_0, \dots,X_{\ell+1})$ partitions $V(T)$.
            \label{claim:structure_nice_decomposition:enum:4}
        \end{enumerate}
    \end{claimlabeled}
    \begin{proofclaim}
        Items~\ref{claim:structure_nice_decomposition:enum:1} and~\ref{claim:structure_nice_decomposition:enum:2} follow easily from the definitions.
        Item~\ref{claim:structure_nice_decomposition:enum:3} follows from~\ref{claim:structure_nice_decomposition:enum:2} and our assumption that $T$ is light.

        For~\ref{claim:structure_nice_decomposition:enum:4}, observe that the sets $(V(P),X_0,\dots,X_{\ell+1})$ are pairwise disjoint by definition. To see that $V(P)\cup X_0 \cup \dots \cup X_{\ell+1} = V(T)$, let $y$ be any vertex in $V(T)\setminus V(P)$. If $y$ does not have any in-neighbour in $V(P)$, then $y$ is an in-neighbour of $v_\ell$, and it belongs to $X_{\ell+1}$. Else if it has at least one out-neighbour $v_i$ in $V(P)$ then, by definition of $X_i$, $y$ belongs to $\bigcup_{j\leq i} X_j$. The claim follows.
    \end{proofclaim}

    Informally, from now on, our goal is to show that we can divide $T$ into smaller pieces $Z_1,\dots,Z_k$ that are easy to colour, and to show that these pieces can be coloured (almost) independently. We then conclude with a dynamic programming argument.
    For this, let us define:
    \begin{align*}
        Y_0 &= X_0 \cup {v_0},\\
        Y_{\ell+1} &= X_{\ell+1} \cup \{v_\ell\}, \text{~~and}\\
        Y_i &= X_i \cup \{v_{i-1},v_i\} \text{~~for every $i\in [1,\ell]$.}
    \end{align*}
    Now each set $Z_j$ should correspond to the union of $7$ consecutive sets $Y_i$s, unless the number of such sets is too small. Formally, let $k=\max(0,\ell-5)$. If $\ell<5$, let 
    \[
        Z_0 = Z_k = \bigcup_{i=0}^{\ell+1} Y_i. 
    \]
    Else if $\ell\geq 5$, for every $i\in [0,k]$, let 
    \[
        Z_i = \bigcup_{j=i}^{i+6} Y_i.
    \]
    In both cases, observe that 
    \[
    \bigcup_{i=0}^{k}Z_i = \bigcup_{i=0}^{\ell+1}Y_i = V(P) \cup \bigcup_{i=0}^{\ell+1}X_i = V(T),
    \]
    where the last equality follows from Claim~\ref{claim:structure_nice_decomposition}\ref{claim:structure_nice_decomposition:enum:4}. 
    Given an index $i\in [0,k]$, a $2$-colouring $\gamma \colon Z_i \to \{1,2\}$ of $T[Z_i]$ is {\bf $i$-coherent} if:
    \begin{enumerate}
        \item $\gamma$ is an acyclic dicolouring of $T[Z_i]$;
        \item for every $u\in Z_i\cap \{s,t\}$, $\phi(u) = 1$; and
        \item for every $u \in Z_i \cap (N^+(s) \cup N^-(t))$, $\phi(u)=2$.
    \end{enumerate}

    Observe that, for any $(s,t)$-colouring $\phi$, its restriction to any set $Z_i$ is $i$-coherent. The first key point is that the converse is also true, hence implying that we only have to find a colouring of $V(T)$ whose restriction to $Z_i$ is $i$-coherent for every $i\in [0, k]$.

    \begin{claimlabeled}[thm:polytime_tour_2colours]
        \label{claim:reduce_to_coherent_colourings}
        A colouring $\phi$ of $T$ is an $(s,t)$-colouring if and only if, for every $i\in [0,k]$, the restriction of $\phi$ to $Z_i$ is $i$-coherent.
    \end{claimlabeled}
    \begin{proofclaim}
        The ``only if'' direction is straightforward. 
        For the opposite direction, assume for a contradiction that $\phi$ is a colouring of $T$ such that, for every $i\in [0,k]$, the restriction $\phi_i$ of $\phi$ to $Z_i$ is $i$-coherent, but $\phi$ is not an $(s,t)$-colouring. 
        
        Since $V(T)=\bigcup_{i=1}^{k}Z_i$, in particular each vertex $u \in \{s,t\} \cup N^+(s) \cup N^-(t)$ belongs to some set $Z_i$. Hence, $\phi_i$ being $i$-coherent, we have $\phi(s) = \phi(t) = 1$ and, for every $u \in N^+(s) \cup N^-(t)$, $\phi(u) = 2$.
        Therefore, since $\phi$ is not an $(s,t)$-colouring, $T$, coloured with $\phi$, contains a directed cycle which is either monochromatic  or alternating between two colour classes. Among all such cycles, we let $\Cscr$ be one of minimum length. Since $T$ is a tournament, $\Cscr$ has length either $3$ (if $\Cscr$ is monochromatic) or $4$ (if $\Cscr$ is alternating).

        Now observe that $V(\Cscr)$ is not included in any set $Z_i$, for otherwise $\phi_i$ would not be an acyclic dicolouring of $T[Z_i]$, a contradiction to $\phi_i$ being $i$-coherent. Therefore, there exist two indices $i,j$ and two vertices $x\in V(\Cscr) \cap Y_i$ and $y\in V(\Cscr)\cap Y_j$ such that $j \geq i+7$. We assume that $x\in X_i\cup \{v_i\}$, for otherwise $x=v_{i-1}$ and $x$ belongs to $Y_{i-1}$, meaning that we can let $i'=i-1$ play the role of $i$. Similarly, we assume that $y\in X_j\cup \{v_{j-1}\}$, for otherwise $y=v_{j}$ and $j'=j+1$ may play the role of $j$.
        
        Since $\Cscr$ has length at most $4$, there exists in $\Cscr$ (and thus in $T$) a directed path $Q$ from $x$ to $y$ of length at most $3$.
        By Claims~\ref{claim:structure_nice_decomposition}\ref{claim:structure_nice_decomposition:enum:1} and~\ref{claim:structure_nice_decomposition:enum:2}, since $x\in X_i \cup \{v_i\}$ and $y\in X_j \cup \{v_{j-1}\}$, we have $x \in \{v_i\} \cup N^+(v_i)$ and $y \in \{v_{j-1}\} \cup N^-(v_{j-1})$. Hence, by adding at most two arcs, $Q$ can be extended into a directed path $Q'$ from $v_i$ to $v_{j-1}$ of length at most $5$. Since the distance from $v_i$ to $v_{j-1}$ along $P$ is at least $j-1-i \geq 6$, this is a contradiction to the choice of $P$.
        The claim follows.
    \end{proofclaim}

    The second key point is that, for every $i\in [0,k]$, the number of $i$-coherent colourings of $T[Z_i]$ is limited.  

    \begin{claimlabeled}[thm:polytime_tour_2colours]
        \label{claim:bound_number_coherent_colourings}
        For every $i\in [0,k]$, there is at most $2^8 \cdot (n+1)^{7}$ distinct $i$-coherent colourings of $T[Z_i]$. Furthermore, the set $\Ccal_i$ of such colourings can be computed in polynomial time.
    \end{claimlabeled}
    \begin{proofclaim}
        Let us fix $i\in [0,k]$. Note that $Z_i$ contains at most $8$ vertices in $V(P)$. Let us fix $\gamma$, one of the $2^8$ possible colourings of $T[Z_i \cap V(P)]$, and evaluate the number of possible extensions of $\gamma$ to $Z_i$. 
        Note that $Z_i$ contains $X_j$ for at most $7$ different indices $j$. We will show that, for each of them, at most $n+1$ colourings of $T[X_j]$ may extend $\gamma$, hence implying the result. For the second part of the statement, it is straightforward to check from the proof below that, in polynomial time, one can iterate over the $2^7 \cdot (n+1)^{6}$ possible colourings, and for each of them check that it is actually $i$-coherent.
        The structure of possible $i$-coherent colourings of $T[Z_i]$ is illustrated in Figure~\ref{fig:i_coherent_colouring}.
        
        \begin{figure}[ht]
            \begin{center}      
                  \begin{tikzpicture}
                    \def\L{1.8};
                    \def\H{1.2};
                    \node[svertex, g-blue, label=left:{$v_{0}$}] (v0) at (0,0) {};
                    \node[svertex, orange, label=below:$v_{1}$] (v1) at (\L,0) {};
                    \node[svertex, orange, label=below:$v_{2}$] (v2) at (2*\L,0) {};
                    \node[svertex, g-blue, label=below:$v_{3}$] (v3) at (3*\L,0) {};
                    \node[svertex, g-blue, label=below:$v_{4}$] (v4) at (4*\L,0) {};
                    \node[svertex, orange, label=below:$v_{5}$] (v5) at (5*\L,0) {};
                    \node[svertex, g-blue, label=below:$v_{6}$] (v6) at (6*\L,0) {};
    
                    \node[Xrect, label=left:$X_0$, fill=orange!25] (X0) at (0,-\H) {};
                    \node[Xrect, label=above:$X_1$] (X1) at (0.5*\L,\H) {};
                    
                    \draw [->,decorate,decoration=snake, thick] (0.5*\L+0.45,\H) -- (0.5*\L-0.46,\H);
                    \node[Xrect, label=above:$X_2$, fill=g-blue!25] (X2) at (1.5*\L,\H) {};
                    \draw [->,decorate,decoration=snake, thick] (1.5*\L+0.45,\H) -- (1.5*\L-0.46,\H);
                    \node[Xrect, label=above:$X_3$] (X3) at (2.5*\L,\H) {};
                    \draw [->,decorate,decoration=snake, thick] (2.5*\L+0.45,\H) -- (2.5*\L-0.46,\H);
                    \node[Xrect, label=above:$X_4$, fill=orange!25] (X4) at (3.5*\L,\H) {};
                    \draw [->,decorate,decoration=snake, thick] (3.5*\L+0.45,\H) -- (3.5*\L-0.46,\H);
                    \node[Xrect, label=above:$X_5$] (X5) at (4.5*\L,\H) {};
                    \draw [->,decorate,decoration=snake, thick] (4.5*\L+0.45,\H) -- (4.5*\L-0.46,\H);
                    \node[Xrect, label=above:$X_6$] (X6) at (5.5*\L,\H) {};
                    \draw [->,decorate,decoration=snake, thick] (5.5*\L+0.45,\H) -- (5.5*\L-0.46,\H);

                    \foreach \i in {0,...,6}{
                        \node[fit= (X\i), inner sep=3pt] (Y\i) {};
                    }

                    \begin{scope}[on background layer]
                        \node[fit=(X1.north west) (X1.-100), fill=g-blue!25, inner sep=0] {};
                        \node[fit=(X1.north east) (X1.-100), fill=orange!25, inner sep=0] {};
                        
                        \node[fit=(X3.north west) (X3.-120), fill=orange!25, inner sep=0] {};
                        \node[fit=(X3.north east) (X3.-120), fill=g-blue!25, inner sep=0] {};
                        
                        \node[fit=(X5.north west) (X5.-80), fill=g-blue!25, inner sep=0] {};
                        \node[fit=(X5.north east) (X5.-80), fill=orange!25, inner sep=0] {};
                        
                        \node[fit=(X6.north west) (X6.-55), fill=orange!25, inner sep=0] {};
                        \node[fit=(X6.north east) (X6.-55), fill=g-blue!25, inner sep=0] {};
                    \end{scope}
                    \foreach \i in {0,1,2,3,4,5,6}{
                        \draw[sarc] (v\i) to (Y\i);
                    }
                    \foreach \i in {0,1,2,3,4,5}{
                        \pgfmathtruncatemacro{\j}{\i+1};
                        \draw[sarc] (v\i) to (v\j);
                    }
                    \draw[sarc, thick, dotted] (v4) to (v5);
                    \foreach \i in {1,2,3,4,5,6}{
                        \pgfmathtruncatemacro{\j}{\i-1};
                        \draw[sarc] (Y\i) to (v\j);
                    }
                  \end{tikzpicture}
              \caption{The only possible structure of any $0$-coherent colouring of $T[Z_0]$. Vertices coloured $1$ are in blue, and vertices coloured $2$ in orange. As indicated by the twisting arrows, the acyclic ordering of each $X_j$, $j\geq 1$, goes from right to left.}
              \label{fig:i_coherent_colouring}
            \end{center}
        \end{figure}
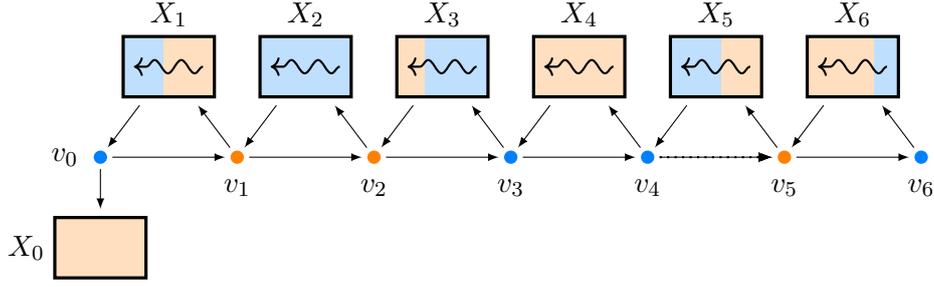
        
        Let us thus fix $0\leq j \leq \ell+1$ such that $X_j \subseteq Z_i$. If $j\in \{0,\ell+1\}$ then, in every $i$-coherent colouring, all vertices in $X_j$ must be coloured $2$, since $X_j \subseteq N^+(s) \cup N^-(t)$.
        Henceforth, assume that $1\leq j \leq \ell$. By Claim~\ref{claim:structure_nice_decomposition}\ref{claim:structure_nice_decomposition:enum:2} and~\ref{claim:structure_nice_decomposition:enum:3}, $X_j \subseteq N^+(v_j) \cap N^-(v_{j-1})$ and $X_j$ is acyclic.

        If $\gamma(v_{j-1}) = \gamma(v_j)=1$, then in every $i$-coherent colouring extending $\gamma$, every vertex $x\in X_j$ must be coloured $2$, for otherwise $(v_{j-1},v_j,x)$ is a monochromatic directed triangle in $T[Z_i]$. Similarly, if $\gamma(v_{j-1}) = \gamma(v_j)=2$, there is only one possible extension of $\gamma$ to $X_j$.

        Assume finally that $1= \gamma(v_{j-1}) \neq \gamma(v_j)=2$, the other remaining case being symmetric. Let $x_1,\dots, x_r$ be the acyclic ordering of $T[X_j]$, where $r=|X_j|$. In every $i$-coherent extension $\beta$ of $\gamma$ to $X_j$, and for every $1\leq a< b\leq r$, if $\beta(x_a) = 1$ then $\beta(x_b)=1$, for otherwise $(v_{j-1},v_j,x_a,x_b)$ is an alternating directed cycle.
        Therefore, $\beta$ is determined by the index $a\in[0,r]$ such that $\phi(x_b)= 1$ for all $b> a$ and $\phi(x_b)=2$ for all $b\leq a$. This shows that at most $r+1\leq n+1$ colourings of $T[X_j]$ may extend $\gamma$, and the claim follows.
    \end{proofclaim}

    We now compute, for every $i\in [k]$, the set $\Ccal_i$ of $i$-coherent colourings of $T[Z_i]$, which can be done in polynomial time by the second part of Claim~\ref{claim:bound_number_coherent_colourings}.
    By Claim~\ref{claim:reduce_to_coherent_colourings}, we know that $(s,t)$ is a nice pair if and only if there exists a colouring $\phi$ of $T$ such that, for every $i\in [0,k]$, the restriction of $\phi$ to $Z_i$ belongs to $\Ccal_i$.

    For every $i \in [0, k]$ and colouring $\phi$ of $T[Z_i]$, we say that $\phi$ is {\bf $(\leq i)$-coherent} if there exists an extension $\phi^\star$ of $\phi$ to $T[\bigcup_{j=0}^{i}Z_j]$ whose restriction to $Z_j$ belongs to $\Ccal_j$ for every $j\in [0,i]$. We further denote by $\Ccal_i^\star$ the set of $(\leq i)$-coherent colourings of $T[Z_i]$. To decide whether $(s,t)$ is a nice pair, it only remains to decide whether $\Ccal_k^\star \neq \emptyset$.

    Observe that, by definition, $\Ccal_0^\star = \Ccal_0$. We can recursively build $\Ccal_i^\star$ from $\Ccal_{i-1}^\star$ as follows.
    \begin{claimlabeled}[thm:polytime_tour_2colours]
        \label{claim:computing_Cistar}
        For every $i\in [1,k]$,
        \[
            \Ccal_i^\star = \{ \phi \in \Ccal_i \mid \mbox{there exists } \psi \in \Ccal_{i-1}^\star \mbox{ such that } \phi(u) = \psi(u)  \text{ for all } u\in Z_i \cap Z_{i-1} \}.
        \]
    \end{claimlabeled}
    \begin{proofclaim}
        Let first $\phi$ be any colouring in $\Ccal_i^\star$, and let $\phi^\star$ be its extension to $T[\bigcup_{j=0}^iZ_j]$. Since $\phi$ is exactly the restriction of $\phi^\star$ to $Z_i$, then $\phi \in \Ccal_i$. Moreover, the restriction $\phi_{i-1}$ of $\phi^\star$ to $Z_{i-1}$ is $(\leq i-1)$-coherent, as witnessed by the restriction of $\phi^\star$ to $\bigcup_{j=0}^{i-1}Z_j$. This shows 
        \[
            \Ccal_i^\star \subseteq \{ \phi \in \Ccal_i \mid \mbox{there exists } \psi \in \Ccal_{i-1}^\star \mbox{ such that } \phi(u) = \psi(u)  \text{ for all } u\in Z_i \cap Z_{i-1} \}.
        \]
        For the opposite direction, let $\phi$ be any $i$-coherent colouring of $T[Z_i]$ such that there exists $\psi \in \Ccal_{i-1}^\star$ agreeing with $\phi$ on $Z_i \cap Z_{i-1}$. Since $\psi$ is $(\leq i-1)$-coherent, there exists an extension $\psi^\star$ of $\psi$ to $\bigcup_{j=0}^{i-1}Z_j$ whose restriction to any set $Z_j$, $j\leq i-1$, is $j$-coherent.

        Let $\phi^\star$ be the colouring of $T[\bigcup_{j=0}^{i}]$ equal to $\phi$ on $Z_i$ and to $\psi^\star$ on $\bigcup_{j=0}^{i-1}Z_j$. This is possible because $\psi^\star(u) = \psi(u) = \phi(u)$ for every $u\in Z_i \cap Z_{i-1}$, and because 
        \begin{equation}
            \label{eq:Zi_cap_Zj}
            Z_i \cap \Big( \bigcup_{j=0}^{i-1}Z_j\Big) = Z_i \cap Z_{i-1}.
        \end{equation}
        To see that \eqref{eq:Zi_cap_Zj} holds, we use the definition of the sets $Z_j$s and the fact that $V(P), X_0,\dots,X_{\ell+1}$ are pairwise disjoint (by Claim~\ref{claim:structure_nice_decomposition}\ref{claim:structure_nice_decomposition:enum:4}). Using these facts, we conclude that:
        \begin{align*}
            Z_i \cap Z_{i-1} &\subseteq Z_i \cap \Big( \bigcup_{j=0}^{i-1}Z_j\Big)\\
            &= \Big(\{ v_{i-1},\dots,v_{i+6}\} \cup \bigcup_{j=i}^{i+6} X_j \Big) \cap \Big(\{ v_{0},\dots,v_{i+5}\} \cup \bigcup_{j=0}^{i+5} X_j  \Big) \\
            &= \{ v_{i-1},\dots,v_{i+5}\} \cup \bigcup_{j=i}^{i+5} X_j \\
            &\subseteq  Z_i \cap Z_{i-1},
        \end{align*}
        implying \eqref{eq:Zi_cap_Zj}.
        It follows by construction that the restriction of $\phi^\star$ to any set $Z_j$, $j\leq i$, is $j$-coherent, and thus $\phi \in \Ccal_i^\star$. The claim follows.
    \end{proofclaim}
    
    It follows from Claim~\ref{claim:computing_Cistar} that, given $\Ccal_{i-1}^\star$, we can build $\Ccal_i$ in polynomial time (since both $\Ccal_i$ and $\Ccal_{i-1}^\star$ have size at most $O(n^6)$ by Claim~\ref{claim:bound_number_coherent_colourings}). After $k\leq n$ steps, in particular, we know whether $\Ccal_k^\star$ is empty, and thus whether $(s,t)$ is a nice pair.
    Moreover, if $(s,t)$ is a nice pair, an $(s,t)$-colouring $\phi$ can be obtained with a simple backtracking, starting from any colouring of $T[Z_{k}]$ in $\Ccal_{k}^\star$ and, for $i=k-1$ to $0$, extending it to $Z_i$ by choosing any colouring in $\Ccal_i^\star$ equal to $\phi$ on $Z_{i+1}$.
\end{proof}

\section{Tight bounds on the acyclic dichromatic number of tournaments}\label{sec:tightbounds}

Recall that, by Theorem~\ref{thm:erdos_moser}, every tournament of order $n$ contains a transitive subtournament of order $\lceil \log n\rceil$. It follows that every tournament of order $n$ admits an acyclic $(n-\lceil \log n\rceil +1)$-dicolouring, which consists of colouring a largest transitive subtournament with one colour, and using one dedicated colour for every other vertex. 
We slightly improve on this trivial upper bound by showing the following.

\begin{theorem}
    \label{thm:upper_bound_T_n}
    Let $T$ be any tournament of order $n$, then 
    \[
        \dict(T) \leq n - \alpha \cdot \log(n) + 27,
    \]
    where $\alpha = \log^{-1}(4/3) > 2.4$.
\end{theorem}

It turns out that, up to the multiplicative constant $\alpha$, Theorem~\ref{thm:upper_bound_T_n} is best possible.

\begin{restatable}{theorem}{thmlowerboundTn}
    \label{thm:lower_bound_T_n}
    There exists $n_0\in \mathbb{N}$ such that, for every integer $n\geq n_0$, there exists a tournament $T$ of order $n$ with
    \[
        \dict(T) \geq n - \beta \cdot \log(n) - \log \log(n),
    \]
    where $\beta = \frac{8}{\log(8/7)} < 60$.
\end{restatable}

Theorems~\ref{thm:lower_bound_T_n} and~\ref{thm:upper_bound_T_n} together raise the following problem.

\begin{problem}
    Find the infimum $\gamma$ of all real numbers $c$ such that there exist infinitely many tournaments $T$ with 
    $
        \dict(T) \geq n- c \cdot \log(n).
    $
    So far, we have $2.4 < \gamma < 60$.
\end{problem}

The remaining of this section is devoted to the proofs of Theorems~\ref{thm:upper_bound_T_n} and~\ref{thm:lower_bound_T_n}.

Our proof of Theorem~\ref{thm:upper_bound_T_n} relies on the notion of {\bf acyclic matchings}. A {\bf matching} $\mathcal{M}$ in a tournament $T$ is a set of pairwise disjoint arcs of $T$. It is {\bf acyclic} if, for every pair of arcs $uv,xy\in \mathcal{M}$, the oriented bipartite graph $T[\{u,v,x,y\}] \setminus \{uv,xy\}$ is acyclic. 

\begin{proof}[Proof of Theorem~\ref{thm:upper_bound_T_n}]
    Let $f\colon x \mapsto \alpha \cdot \log(x) - 27(1- \frac{1}{x})$ be defined over $\mathbb{R}^+$.
    We prove by induction on $n\in \mathbb{N}^\ast$ that every tournament $T$ on $n$ vertices admits an acyclic matching $\mathcal{M}$ of size at least $f(n)$.
    The statement then follows, as an acyclic $(n-|\mathcal{M}|)$-dicolouring of $T$ can be obtained by giving a common colour to any pair of vertices matched in $\mathcal{M}$, and   a private colour to  each of the remaining vertices.

    Observe that $f(n)\leq 0$ when $n\leq 11$, so the result trivially holds in this case. We thus assume that $n\geq 12$. 
    Let $uv \in A(T)$ be an arc such that either $u$ is a source or 
    \[
        |N^+(u) \cap N^-(v)| \leq \frac{n-3}{4}.
    \]
    To see that such an arc exists, let $v$ be any vertex of $T$ with minimum in-degree, and assume that $d^-(v) \geq 1$ (otherwise $vw$ may play the role of $uv$ for any $w\in N^+(v)$). Since $\sum_{v\in V(T)}d^-(v) = \binom{n}{2}$, $d^-(v) \leq \frac{n-1}{2}$.
    Let $u\in N^-(v)$ be any vertex whose out-degree is minimum in $T[N^-(v)]$. Again, $u$ has out-degree at most $\frac{1}{2}(\frac{n-1}{2}-1)$ in $T[N^-(v)]$, and it follows that, in $T$, 
    \[
        |N^+(u) \cap N^-(v)| \leq \frac{n-3}{4}.
    \]

    If $u$ is a source, let $X= \emptyset$, and let $X=N^+(u) \cap N^-(v)$ otherwise. Let  $T'$ be the tournament $T-(X\cup \{u,v\})$. Note that $|V(T')| \geq \lceil \frac{1}{4}(3n-5) \rceil \geq 8$, and that $f$ is increasing over $[8,+\infty)$. Therefore, by induction, $T'$ admits an acyclic matching $\mathcal{M}'$ of size at least $f\left(\frac{3n-5}{4}\right)$.
    
    Let $\mathcal{M} = \mathcal{M}' \cup \{uv\}$. We claim that $\mathcal{M}$ is acyclic. To see this, let $xy$ be any arc of $\mathcal{M}'$, we only have to show that $C = T[\{u,v,x,y\}] \setminus \{xy,uv\}$ is acyclic. Note that $C$ is an orientation of $C_4$, so it is acyclic if and only if it contains a source or a sink. Therefore, if $x\in N^+(u) \cap N^+(v)$ or $x\in N^-(u) \cap N^-(v)$, then $C$ is acyclic. By construction, of $T'$, we thus have $x\in N^-(u) \cap N^+(v)$ and analogously $y\in N^-(u) \cap N^+(v)$. Hence, $u$ is a sink and $C$ is acyclic, as desired.

    It remains to check that $1+f\left(\frac{3n-5}{4}\right) \geq f(n)$. Note first that, for any real $x>5/4$, we have
    \begin{align*}
        \log(x) - \log\left(x-\frac{5}{4}\right) = \log\left(1+ \frac{5}{4x-5}\right) \leq \frac{5}{(4x-5)\ln(2)},
    \end{align*}
    where in the last inequality we used that $\ln(1+x) \leq x$ and that $\log(1+x) = \frac{1}{\ln(2)}\ln(1+x)$. 
    It follows that 
    \begin{align*}
        |\mathcal{M}| \geq 1+ f\left(\frac{3n-5}{4}\right) &=1+\alpha \cdot  \log\left(\frac{3n-5}{4}\right) - 27\left(1-\frac{4}{3n-5}\right)\\
        &\geq 1+\alpha \cdot \log\left(\frac{3}{4}n\right) - \frac{5\alpha}{(3n-5)\ln(2)} - 27\left(1-\frac{4}{3n-5}\right)\\
        &= \alpha \cdot \log(n) +  \underbrace{\Big(1- \alpha \cdot \log\left(4/3\right)\Big)}_{=0} - \frac{5\alpha}{(3n-5)\ln(2)} - 27\left(1-\frac{4}{3n-5}\right)\\
        &\geq \alpha \cdot \log(n) - \frac{9}{n} - 27 + \frac{36}{n}\\
        &= f(n).
    \end{align*}
    The result follows.
\end{proof}

We now prove \Cref{thm:lower_bound_T_n}. Our proof relies on probabilistic arguments. To minimize the value of $\beta$ we obtain, we first evaluate the probabilities that random orientations of specific small graphs are acyclic. For this, we make use of the following well-known result of Stanley~\cite{stanleyDM5}.
For every graph $G$, there exists a unique polynomial $\Pscr(G,x)$, called the {\bf chromatic polynomial} of $G$, with the property that, for every $k\in \mathbb{N}$, the number of proper $k$-colourings of $G$ is exactly $\Pscr(G,k)$, see for instance~\cite[Section~14.7]{bondy1976}. Stanley's theorem~\cite{stanleyDM5} draws a connection between the chromatic polynomial of a graph and the number of its acyclic orientations.

\begin{theorem}[Stanley~\cite{stanleyDM5}]
    \label{thm:stanley}
    Every graph $G$ admits $|\Pscr(G,-1)|$ distinct acyclic orientations.
\end{theorem}

We derive the following. 

\begin{lemma}
    \label{lemma:prob_acyclic}
    There exist 14 acyclic orientations of $K_{2,2}$, 46 acyclic orientations of $K_{2,3}$, and 230 acyclic orientations of $K_{3,3}$. In particular, if $D_{2,2}$, $D_{2,3}$, and $D_{3,3}$ are random orientations of $K_{2,2}$, $K_{2,3}$, and $K_{3,3}$ respectively, then 
    \[  
        \PP(D_{2,2} \mbox{\rm ~is acyclic}) = \frac{7}{8} \mbox{,~~~}
        \PP(D_{2,3} \mbox{\rm ~is acyclic}) = \frac{23}{32} \mbox{,~~~and~~~} \PP(D_{3,3} \mbox{\rm ~is acyclic}) = \frac{115}{256}.
    \]
\end{lemma}
\begin{proof}
    Straightforward counting arguments imply:
    \begin{align*}
        \Pscr(K_{2,2},x)&= x^4-4x^3+6x^2-3x,\\
        \Pscr(K_{2,3},x) &= x^5 - 6 x^4 + 15 x^3 - 17 x^2 + 7 x, \\
        \text{and~~} \Pscr(K_{3,3},x) &= x^6 - 9x^5 + 36x^4 - 75x^3 + 78x^2 - 31x.
    \end{align*}
    By~\Cref{thm:stanley}, for integers $n,m$, the number $\Ascr_{n,m}$ of acyclic orientations of $K_{n,m}$ is precisely $|\Pscr(K_{n,m},-1)|$.
    It follows that 
    \[
        \Ascr_{2,2} (x) = 14, \text{~~}
        \Ascr_{2,3} (x) = 46,
        \text{~~and~~} \Ascr_{3,3} (x) = 230,
    \]
    as desired. The second part of the statement follows from the observation that $K_{n,m}$ admits $2^{nm}$ distinct orientations.
\end{proof}

We slightly extend the notion of acyclic matching as follows. An {\bf acyclic partition} in a tournament $T$ is a collection $(X_1,\dots,X_t)$ of subsets of $V(T)$ such that:
\begin{enumerate}[label=$(\roman*)$]
    \item for each $i\in [t]$, $2\leq |X_i|\leq 3$;
    \item the sets $X_1,\dots, X_t$ are pairwise disjoint; and
    \item for all $1\leq i<j\leq t$, the oriented bipartite graph $T[X_i,X_j]$ is acyclic. 
\end{enumerate}

We say that a subset $X\subseteq V(T)$ of vertices is {\bf partitionable} if there exists an acyclic partition $(X_1,\dots,X_t)$ such that $X = \bigcup_{i=1}^tX_i$.
We now prove Theorem~\ref{thm:lower_bound_T_n}, which we recall here for convenience.

\thmlowerboundTn*
\begin{proof}
    We do not give the exact value of $n_0$, we simply assume that it is large enough so that all upcoming inequalities hold.

    Let $T$ be a random tournament on $n$ vertices (that is, each edge $\{u,v\}$ is oriented independently and  uniformly at random). Let $k = \beta  \cdot \log(n) + \log\log(n)$. We first show that, with positive probability, none of the $k$-sets $X\subseteq V(T)$ is partitionable. We then conclude that, with positive probability, $\dict(D) \geq n-k$, hence implying the result. 

    Let us fix a pairwise disjoint collection $(X_1,\dots,X_t)$ of sets of size either $2$ or $3$, and whose union has size $\ell =\lceil k \rceil$. We first bound the probability that $(X_1,\dots,X_t)$ is an acyclic partition in $T$. For any $1\leq i<j \leq t$, let $E_{i,j}$ denote the event that the bipartite oriented graph $T[X_i,X_j]$ is acyclic. By \Cref{lemma:prob_acyclic}, we have 
    \begin{itemize}
        \item $\PP(E_{i,j}) = \frac{7}{8}$ if $|X_i| = |X_j| = 2$,
        \item $\PP(E_{i,j}) =\frac{23}{32}$ if $|X_i| \neq |X_j|$, and
        \item $\PP(E_{i,j}) =\frac{115}{256}$ if $|X_i| = |X_j|=3$.
    \end{itemize}
    For any distinct pairs $\{i,j\} \subseteq [k]$ and $\{r,s\} \subseteq [k]$, observe that the events $E_{i,j}$ and $E_{r,\ell}$ are independent. 
    It follows that 
    \[
        \PP\Big((X_1,\dots,X_t) \mbox{ is an acyclic partition}\Big) \leq 
        \underbrace{\left(\frac{7}{8}\right)^{\frac{1}{2}k_2(k_2-1)} \cdot \left(\frac{23}{32}\right)^{k_2 k_3} \cdot \left(\frac{115}{256}\right)^{\frac{1}{2}k_3(k_3-1)}}_{g(k_2,k_3)=},
    \]
    where $k_2$, respectively $k_3$, denotes the number of sets $X_i$ of size 2, respectively 3 (So $\ell=2k_2+3k_3$).
    Let us simplify the upper bound above. Note that $g(k_2,k_3)$ denotes the right-hand side above.
    Let us fix two arbitrary real numbers $\ell_2,\ell_3\geq 0$ such that $2\ell_2+3\ell_3 = \ell$. When $\ell$ is large enough, we have
    \begin{align*}
        &g(\ell_2+3,\ell_3-2)\\
        &=  \left(\frac{7}{8}\right)^{\frac{1}{2}(\ell_2+3)(\ell_2+2)} 
            \cdot \left(\frac{23}{32}\right)^{(\ell_2+3)(\ell_3-2)} \cdot \left(\frac{115}{256}\right)^{\frac{1}{2}(\ell_3-2)(\ell_3-3)}\\
        &= g(\ell_2,\ell_3) \cdot \left(\frac{7}{8}\right)^{3\ell_2+3} \cdot \left(\frac{23}{32}\right)^{3(\ell_3-2)-2\ell_2} \cdot \left(\frac{115}{256}\right)^{-2\ell_3+3}\\
        &= g(\ell_2,\ell_3)\cdot \left(\underbrace{\left(\frac{7}{8}\right)^{3} \left(\frac{23}{32}\right)^{-2} }_{\approx 1.3}\right)^{\ell_2}
        \cdot \left(\underbrace{\left(\frac{23}{32}\right)^{3} \left(\frac{115}{256}\right)^{-2}}_{\approx 1.84} \right)^{\ell_3}
        \cdot \underbrace{\left(\frac{7}{8}\right)^{3} \left(\frac{23}{32}\right)^{-6} \left(\frac{115}{256}\right)^{3}}_{\approx 0.44} \\
        &>g(\ell_2,\ell_3),
    \end{align*}
    where in the last inequality we used that either $\ell_2$ or $\ell_3$ is large enough.
    Since $2k_2+3k_3=\ell$, it follows from the remark above and a straightforward recurrence that either $\ell$ is even and $g(k_2,k_3) \leq g(\ell/2,0)$ or $\ell$ is odd and $g(k_2,k_3) \leq g(\frac{\ell-3}{2},1)$. 
    Note that 
    \[
        g(\ell/2,0) = \left(\frac{7}{8}\right)^{\frac{1}{8}\ell(\ell-2)}<  \left(\frac{7}{8}\right)^{\frac{1}{8}\ell^2 - 2\ell},
    \]
    and that 
    \[
        g\left(\frac{\ell-3}{2},1\right) < \left(\frac{7}{8}\right)^{\frac{1}{8}(\ell^2-8\ell+15)} < \left(\frac{7}{8}\right)^{\frac{1}{8}\ell^2 - 2\ell},
    \]
    where in the last inequality we used that $n$ (and thus $k$ and $\ell$) is large enough, so $\ell > 15/8$.
    In both cases, we conclude that 
    \[
        \PP\Big((X_1,\dots,X_t) \mbox{ is an acyclic partition}\Big) < \left(\frac{7}{8}\right)^{\frac{1}{8}\ell^2 - 2\ell}  = 2^{-\frac{1}{\beta}\ell^2 +16\beta \ell}.
    \]
    For a fixed set of vertices $X$ of size $\ell$, there are at most $\ell^\ell = 2^{\ell \log (\ell)}$ partitions of $X$. It thus follows from the union bound that 
    \[
        \PP(X \mbox{ is partitionable}) < 2^{-\frac{1}{\beta}\ell^2 + \ell \log(\ell) +16\beta \ell}.
    \]
    Finally, there exist exactly $\binom{n}{\ell} < n^\ell = 2^{\ell \log(n)}$ sets of vertices of size $\ell$. By the union bound again, it follows that
    \[
        \PP\Big(\exists X \subseteq V(T), |X|=\ell, \text{ $X$ is partitionable}\Big) < 2^{\ell\log(n) -\frac{1}{\beta}\ell^2 + \ell \log(\ell) +16\beta \ell}.
    \]
    Recall that $k \leq \ell \leq k+1$. Hence
    \begin{align*}
        &\ell\log(n) -\frac{1}{\beta}\ell^2 + \ell \log(\ell) +16\beta \ell\\
        &\leq k\log(n) - \frac{1}{\beta}k^2 + k \log(k) + O(\log n)\\
        &\leq (1-\beta)\log(n) \log \log(n) + O(\log(n))
    \end{align*}
    Since $\beta>1$ and because $n$ is large enough, it follows that 
    \[
    \PP\Big(\exists X \subseteq V(T), |X|=\ell, \text{ $X$ is partitionable}\Big) < \frac{1}{2}.
    \]
    It is straightforward to check that the exact same argument applied with $\ell+1$ instead of $\ell$ holds, hence implying that 
    \[
    \PP\Big(\exists X \subseteq V(T), |X|=\ell+1, \text{ $X$ is partitionable}\Big) < \frac{1}{2}.
    \]
    By the union bound, we conclude on the existence of a tournament in which none of the sets of size $\ell$ or $\ell+1$ is partitionable. We fix such a tournament $T$, and claim that $\dict(T) \geq n-\ell+1$, hence implying the result. Assume for a contradiction that $\dict(T)\leq n-\ell$, and let $\phi$ be an acyclic dicolouring of $T$ using at most $n-\ell$ colours.

    Let $Y_1,\dots,Y_{r}$ be the colour classes of $\phi$, labelled such that $|Y_i| \geq |Y_{i+1}|$ for every $1\leq i \leq n-\ell-1$. Note that $|Y_1|\geq 2$, otherwise $\phi$ uses $n$ colours. Let $h$ be the largest index such that $|Y_i| \geq 2$ and let $Y= \bigcup_{i=1}^h Y_i$. 
    Observe that the number of colours used by $\phi$ is exactly
    \[
        n-\sum_{i=1}^h (|Y_i|-1) \leq n-\ell,
    \]
    hence implying $|Y| \geq \ell+1$. By successively removing from $Y$ either a single vertex in a colour class of size at least $3$ or both vertices in a colour class of size exactly $2$, we conclude on the existence of $X\subseteq Y$ and a partition $(X_1,\dots,X_{h'})$ of $X$ such that:
    \begin{itemize}
        \item $\ell \leq |X| \leq \ell+1$,
        \item for every $i\in [h']$, $X_i \subseteq Y_i$, and
        \item for every $i\in [h']$, $|X_i|\geq 2$.
    \end{itemize}
    For every $i\in [h']$, we let $(X_i^1,\dots,X_i^{\iota_i})$ be an arbitrary partition of $X_i$ into sets of size either $2$ or $3$ (which necessarily exists as $|X_i| \geq 2$).
    We consider the partition
    \[
    (X_1^1,\dots,X_1^{\iota_1}, X_2^1,\dots,X_{h'-1}^{\iota_{h'-1}}, X_{h'}^1,\dots,X_{h'}^{\iota_{h'}}),
    \]
    that we relabel $Z_1,\dots,Z_{\zeta}$ for convenience. We finally argue that $(Z_i)_{i\in [\zeta]}$ is an acyclic partition, hence showing that $X$ is partitionable, a contradiction as $|X|\in \{\ell,\ell+1\}$.

    To see this, let $1\leq i<j\leq \zeta$. By construction, there exist $1\leq i',j' \leq h$ such that $Z_i\subseteq Y_{i'}$ and $Z_j \subseteq Z_{j'}$. If $i'\neq j'$ then $T[Z_i,Z_j]$ is a subdigraph of $T[Y_{i'},Y_{j'}]$, and it is acyclic because $\phi$ is an acyclic dicolouring. If $i'=j'$, then $T[Z_i,Z_j]$ is a subdigraph of $T[Y_{i'}]$, and  it is acyclic because $\phi$ is an acyclic dicolouring. The result follows.
\end{proof}

\section{Degenerate graphs}\label{sec:degenerate}

Let $d$ be an integer. A graph $G$ is {\bf $d$-degenerate} if every non-empty subgraph $H$ of $G$ contains a vertex of degree at most $d$. 
It is well-known that $2$-degenerate graphs have unbounded acyclic chromatic number (see for instance~\cite[Corollary 3]{woodDMTCS7}), while $1$-degenerate graphs ({\it i.e.} forests) trivially admit an acyclic $2$-colouring. 

In this section, we show a similar behaviour in the directed setting. A digraph is $d$-degenerate if its underlying graph is.
Given a digraph $D$ and a partition $(X_1,X_2)$ of $V(D)$, we denote by $D^{rev}(X_1,X_2)$ the digraph $D$ after reversing the direction of all arcs between $X_1$ and $X_2$.
We will show the following theorem.

\begin{theorem} \label{thm:MainAcyclicI}
For every $2$-degenerate oriented graph $D$, there exists a partition $(X_1,X_2)$ of $V(D)$
such that $D^{rev}(X_1,X_2)$ is acyclic.
\end{theorem}

\begin{proof}
  We will prove the theorem by induction on the order, $n$, of $D$. Clearly the theorem holds when the order is at most $3$, since in this case we can use the partition $(\emptyset,V(D))$ when $D$ is acyclic and $(\{v\},V(D)\setminus \{v\})$ for 
any $v \in V(D)$ when $D$ is not acyclic (and therefore $D$ is $3$-cycle). So assume $n \geq 4$ and the theorem 
holds for all smaller values of $n$.

Since $D$ is $2$-degenerate, there exists a vertex $x \in V(D)$ with $d(x) \leq 2$. Let $D' = D - x$ and let
$(X_1',X_2')$ be a partition of $V(D')$ such that $(D')^{rev}(X_1',X_2')$ is acyclic, which exists by our induction hypothesis.
If $d(x) \leq 1$ we can add $x$ to either $X_1'$ or $X_2'$ in order to get the desired partition of $D$ (as $x$ cannot belong to any cycles), so we may assume that $d(x)=2$. We now consider the following cases, which exhaust all possibilities.

\2

{\bf Case 1.} {\em $d^+(x)=d^-(x)=1$ and both neighbours of $x$ belong to the same $X_i'$.}

In this case, we may without loss of generality assume that $v_1 x v_2$ is a path in $D$ and $v_1,v_2 \in X_1'$.
If there is a $(v_1,v_2)$-path, $P$, in $(D')^{rev}(X_1',X_2')$ then we will let $X_1=X_1' \cup \{x\}$ and $X_2=X_2'$.
If there is a cycle in $D^{rev}(X_1,X_2)$, then the cycle must use the path $v_1 x v_2$. But substituting 
$v_1 x v_2$ in the cycle by $P$ gives us a closed walk in $(D')^{rev}(X_1',X_2')$ and therefore also a cycle
in $(D')^{rev}(X_1',X_2')$, a contradiction. So $D^{rev}(X_1,X_2)$ is acyclic in this case.

So we may assume that there is no $(v_1,v_2)$-path, $P$, in $(D')^{rev}(X_1',X_2')$. In this case, let
$X_1=X_1'$ and $X_2=X_2' \cup \{x\}$.
If there now is a cycle in $D^{rev}(X_1,X_2)$, then the cycle must use the path $v_2 x v_1$ and a path from $v_1$ to $v_2$. 
However, we assumed there was no $(v_1,v_2)$-path in $(D')^{rev}(X_1',X_2')$, a contradiction. So in Case~1 we can always 
find the desired partition.

\2

{\bf Case 2.} {\em $d^+(x)=d^-(x)=1$ and the neighbours of $x$ belong to different sets in $(X_1',X_2')$.}

In this case, let $X_1=X_1' \cup \{x\}$ and $X_2=X_2'$ (or $X_1=X_1'$ and $X_2=X_2' \cup \{x\}$ would also work).
Now the out-degree or in-degree of $x$ in $D^{rev}(X_1,X_2)$ will be zero so  $D^{rev}(X_1,X_2)$ will be acyclic.

\2

{\bf Case 3.} {\em $\{d^+(x),d^-(x)\}=\{0,2\}$ and both neighbours of $x$ belong to the same $X_i'$.}

In this case, let $X_1=X_1' \cup \{x\}$ and $X_2=X_2'$ (or $X_1=X_1'$ and $X_2=X_2' \cup \{x\}$ would also work).
Now the out-degree or in-degree of $x$ in $D^{rev}(X_1,X_2)$ will be zero so  $D^{rev}(X_1,X_2)$ will be acyclic.

\2

{\bf Case 4.} {\em $\{d^+(x),d^-(x)\}=\{0,2\}$ and the neighbours of $x$ belong to different sets in $(X_1',X_2')$.}

The proof of this case follows similar lines to that of Case~1. 
In this case, assume without loss of generality, and by directional duality, that $N^+(x)=\{v_1,v_2\}$ and $v_i \in X_i'$ for $i=1,2$.

If there is a $(v_2,v_1)$-path, $P$, in $(D')^{rev}(X_1',X_2')$ then we will let $X_1=X_1' \cup \{x\}$ and $X_2=X_2'$.
If there is a cycle in $D^{rev}(X_1,X_2)$, then the cycle must use the path $v_2 x v_1$. But substituting 
$v_2 x v_1$ in the cycle by $P$ gives us a closed walk in $(D')^{rev}(X_1',X_2')$ and therefore also a cycle
in $(D')^{rev}(X_1',X_2')$, a contradiction. So $D^{rev}(X_1,X_2)$ is acyclic in this case.

So we may assume that there is no $(v_2,v_1)$-path in $(D')^{rev}(X_1',X_2')$. In this case, let
$X_1=X_1'$ and $X_2=X_2' \cup \{x\}$.
If there now is a cycle in $D^{rev}(X_1,X_2)$, then the cycle must use the path $v_1 x v_2$ and a path from $v_2$ to $v_1$ in $D'$. 
However, we assumed there was no $(v_2,v_1)$-path in $(D')^{rev}(X_1',X_2')$, a contradiction. So in Case~4 we can always 
find the desired partition, and the result follows.
\end{proof}

\begin{corollary} \label{cor:MainAcyclicI} 
Every $2$-degenerate oriented graph satisfies $\dict(D)\leq 2$.
\end{corollary}

\begin{proof}  
  By Theorem~\ref{thm:MainAcyclicI}, we note that
there exists a partition $(X_1,X_2)$ of $V(D)$ such that $D^{rev}(X_1,X_2)$ is acyclic.
However, this implies that $D[X_1]$ and $D[X_2]$ and $D[X_1,X_2]$ must all be acyclic as well.
So $(X_1,X_2)$ is the desired colouring.
\end{proof} 

In particular, we derive the following, as outerplanar graphs are $2$-degenerate (see for instance~\cite[Exercise~10.2.11]{bondy1976}). 

\begin{corollary}
    Every oriented outerplanar graph $D$ satisfies $\dict(D) \leq 2$.
\end{corollary}

We now show that Corollary~\ref{cor:MainAcyclicI} is best possible by constructing $3$-degenerate oriented graphs with arbitrarily large acyclic dichromatic number. We make use of the celebrated result of Ramsey in the following general form, which can be found in many basic textbooks on graph theory, see for instance~\cite[Theorem~9.1.3]{diestel2017}.

\begin{theorem}[{\sc Ramsey}~\cite{ramsey1930}]
    \label{thm:ramsey}
    For all $k,b,r\geq 1$, there exists an integer $n\geq r$ such that, for every mapping $\phi\colon \binom{[n]}{b}\to [k]$, there exists $X\in \binom{[n]}{r}$ such that $\phi(x)=\phi(y)$ for all $x,y\in \binom{X}{b}$.
\end{theorem}

For integers $k,b,r\geq 1$, we denote by $R_k(b,r)$ the smallest integer $n$ for which Theorem~\ref{thm:ramsey} holds.

\begin{theorem}
    For every $k\in \mathbb{N}$, there exists a $3$-degenerate oriented graph $D$ with $\dict(D) \geq k$.
\end{theorem}
\begin{proof}
    Let us fix an arbitrary integer $k\in \mathbb{N}^\ast$ and build a $3$-degenerate oriented graph $D$ with $\dict(D) \geq k$. We let $D$ be the oriented graph obtained as follows. Start with an independent set $S$ of size $(k-1)\cdot R_{k-1}(3,4)$, and take an arbitrary acyclic order $\prec$ on $S$. For every triplet $\{u,v,w\}\subseteq S$ with $u\prec v\prec w$, we add a vertex $x_{u,v,w}$ and the arcs $x_{u,v,w}u$,$x_{u,v,w}w$, and $vx_{u,v,w}$. 

    Clearly, $D$ is $3$-degenerate since every vertex $x\in V(D)\setminus S$ has degree $3$, and $S$ is an independent set. Assume for a contradiction that $\dict(D) \leq k-1$, and let $\phi$ be an acyclic $(k-1)$-dicolouring of $D$. By the pigeonhole principle, there exists a colour $i\in [k-1]$ such that $S'=|S\cap \phi^{-1}(i)| \geq R_{k-1}(3,4)$.

    We define $\psi \colon \binom{S'}{3} \to [k-1]$ as follows: for every triplet $\{u,v,w\}\subseteq S'$ with $u\prec v \prec w$, let $\psi(\{u,v,w\}) = \phi(x_{u,v,w})$. By Theorem~\ref{thm:ramsey} and by definition of $R_{k-1}(3,4)$, there exists $S'' = \{s_1,s_2,s_3,s_4\} \subseteq S'$ such that $\psi$ is constant on $\binom{S'}{3}$. Assume without loss of generality that $s_1 \prec s_2 \prec s_3 \prec s_4$. Then $C=(s_2,x_{s_1,s_2,s_3},s_3,x_{s_2,s_3,s_4},s_2)$ is a directed $4$-cycle of $D$ such that $\phi(s_2)=\phi(s_4)$ (because $s_2,s_4\in S'$) and $\phi(x_{s_1,s_2,s_3}) = \phi(x_{s_2,s_3,s_4}) = \psi(\{s_1,s_2,s_3\})$.
    Therefore, $C$ is either monochromatic or alternating between two colours, a contradiction to the choice of~$\phi$.
\end{proof}

\section{Remarks and further open problems}\label{sec:remarks}

We conclude with a branch of further open problems. We start with the case of planar graphs. In 1979, Borodin~\cite{borodinDM25} obtained the following celebrated result, resolving a conjecture of Grünbaum~\cite{grunbaumIJM14}.

\begin{theorem}[Borodin~\cite{borodinDM25}]
    Every planar graph admits a proper $5$-colouring such that the union of any two colour classes induces a forest.
\end{theorem}

The same bound directly follows for the acyclic dichromatic number.

\begin{corollary}
    \label{cor:planar}
    Every planar oriented graph $D$ satisfies $\dict(D) \leq 5$.
\end{corollary}

\begin{figure}[ht]
    \centering
    \begin{tikzpicture}
    \def\L{1.2}
    \node[labelledvertex, SkyBlue] (u) at (3.5*\L,-3) {$x$};
    \node[labelledvertex, SkyBlue] (v) at (3.5*\L,3) {$y$};
    \foreach \i in {1,2,3}{
        \node[labelledvertex, purple] (x\i) at (\i*\L,0) {\footnotesize $z_\i$};
    }
    \foreach \i in {4,5,6}{
        \node[labelledvertex, purple] (x\i) at (\i*\L,0) {\footnotesize $z_\i$};
    }
    
    \draw[thick] (u) to[out=0, in=-90] (8,0);
    \draw[arc] (8,0) to[out=90, in=0] (v);
    
    \draw[arc] (x2) to (x1);
    \draw[arc] (x2) to (x3);
    \draw[arc] (x4) to (x3);
    \draw[arc] (x5) to (x4);
    \draw[arc] (x5) to (x6);
    
    \foreach \i in {1,3,6}{
        \draw[arc] (x\i) to (u);
        \draw[arc] (x\i) to (v);
    }
    \foreach \i in {2,5}{
        \draw[arc] (u) to (x\i);
        \draw[arc] (v) to (x\i);
    }
    \draw[arc] (u) to (x4);
    \draw[arc] (x4) to (v);
    
    \end{tikzpicture}
    \caption{The digraph $W[x,y]$ together with one of its acyclic $2$-dicolourings.}
    \label{fig:Wuv}
\end{figure}
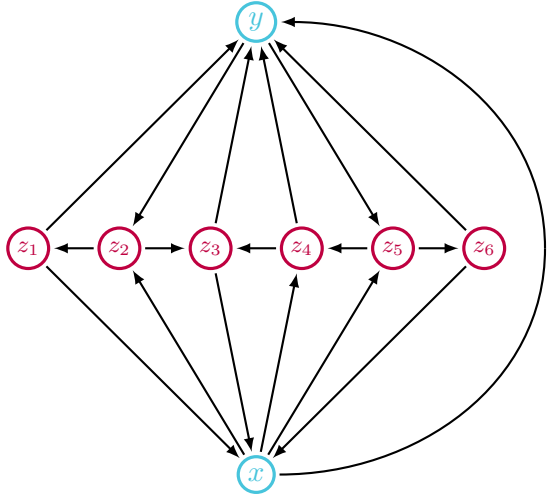

\begin{proposition}
    \label{prop:dict3}
    There exists a planar oriented graph $D$ with $\dict(D) \geq 3$.
\end{proposition}
\begin{proof}
    Consider the digraph $D$ obtained from a directed triangle $(u,v,w)$ by further gluing $W[u,v]$, $W[v,w]$, and $W[w,u]$, where $W[x,y]$ is the planar oriented graph given in Figure~\ref{fig:Wuv}. The statement follows from the observation that, in every $2$-dicolouring of $W[x,y]$, both $x$ and $y$ receive the same colour.
    
    Assume for a contradiction that $\phi$ is an acyclic $2$-dicolouring of $W[x,y]$ with $\phi(x) = 1$ and $\phi(y) = 2$. Then $z_1$ and $z_2$ cannot receive the same colour, for otherwise either $\{z_1,x,z_2\}$ or $\{z_1,y,z_2\}$ induces a monochromatic directed cycle. Moreover, we have $\phi(z_1)=1$ and $\phi(z_2) = 2$, for otherwise $(x,y,z_2,z_1)$ is a directed cycle alternating between colours $1$ and $2$. 
    Hence $\phi(z_3)=1$, otherwise $(y,z_2,z_3)$ is monochromatic.
    By symmetry, we obtain that $\phi(z_6)=1$, $\phi(z_5)=2$, and $\phi(z_4)=1$. It follows that $(x,z_4,z_3)$ is monochromatic, a contradiction. 
\end{proof}

We do not believe that the bound of Corollary~\ref{cor:planar} is tight and pose the following.

\begin{conjecture}
    Every oriented planar graph $D$ satisfies $\dict(D) \leq 3$.
\end{conjecture}

We conclude with several complexity problems. The first one was already mentioned in Section~\ref{subsec:complexity:tournaments}.

\begin{conjecture}
    For every fixed $k\geq 3$ it is polynomial to decide whether a tournament $T$ has $\dict(T)\leq k$
\end{conjecture}

We proved in Section~\ref{sec:complexity} that there exists a polynomial-time algorithm deciding whether a tournament has an acyclic $2$-dicolouring. We wonder if such an algorithm can be extended to cover larger classes of oriented graph.

\begin{problem}
\label{prob:smallalpha}
What is the complexity, for every fixed $k\in \mathbb{N}$, of deciding if a digraph $D$ with $\alpha(D)\leq k$ has $\dict(D)\leq 2$? 
\end{problem}

Already the case when $\alpha(D)=2$ is open, even for digraphs whose vertex set can be covered by two semicomplete digraphs.

A digraph $D$ is a {\bf locally tournament digraph} if $D[N^+[v]]$ and $D[N^-[v]]$ are tournaments for every vertex $v$ of $D$. These digraphs, introduced in \cite{bangJGT14}, form a rich generalization of tournaments and share many properties with these \cite[Chapter 6]{bang2018}.

\begin{problem}
    What is the complexity of deciding $\dict(D)\leq 2$ when $D$ is a locally tournament digraph?
\end{problem}

Note that the definition of light still works for locally tournament digraphs and  the proof of Claim~\ref{claim:T_light} carries over to locally tournament digraphs.

Another generalization of tournaments is the class of {\bf multipartite tournaments} which are orientations of complete multipartite graphs. 
\begin{problem}
    What is the complexity of deciding whether a multipartite tournament has acyclic dichromatic number at most 2?
\end{problem}

This is open already for bipartite tournaments, but polynomial for the special class of multipartite tournaments called {\bf extended tournaments}. These oriented graphs are of the form $D=T[I_{k_1},I_{k_2},\ldots{},I_{k_t}]$ where $T$ is a tournament and $D$ is obtained from $T$  by substituting an independent set of vertices $I_{k_j}$ for the $j$th vertex of $T$ for $j\in [|V(T)|]$. It is easy to check that if $D$ is an extension of $T$ then $\dict(D)=\dict(T)$.

Another generalization of tournaments is the class of {\bf quasi-transitive digraphs}. These are digraphs containing no induced directed path of length 2. They have a very useful recursive structure.

\begin{theorem}\cite{bangJGT20}
  \label{thm:qtdchar}
  Let $D$ be a quasi-transitive digraph.
  \begin{itemize}
  \item[(a)] If $D$ is not strongly connected, then there exists a transitive oriented graph $T$ such that $D=T[H_1,H_2,\ldots{},H_t]$, where each $H_i$ is a strongly connected quasi-transitive digraph.
  \item[(b)] If $D$ is strongly connected, then there exists a strong semicomplete digraph $S$ such that
    $D=S[Q_1,Q_2,\ldots{},Q_s]$ where each $Q_i$ is either a single vertex or a non-strong quasi-transitive digraph. 
  \end{itemize}
  \end{theorem}

  Note that a strong quasi-transitive digraph $D=S[Q_1,Q_2,\ldots{},Q_s]$ contains the semicomplete digraph $S$ as an induced subdigraph (just take an arbitrary vertex from each $Q_i$). Thus, $\dict(D)\leq 2$ only if we also have $\dict(S)\leq 2$. The other direction does not always hold, e.g. in the case when each $Q_i$ contains a directed cycle. In fact, the following holds.
  
  \begin{theorem}
  A strong quasi-transitive digraph $D=S[Q_1,Q_2,\ldots{},Q_s]$ 
  has $\dict(D)=2$ if and only if $\dict(S)=2$ and each $Q_i$ is acyclic. 
  \end{theorem}
  \begin{proof}
      Let $D=S[Q_1,Q_2,\ldots{},Q_s]$ be a strong quasi-transitive digraph without 2-cycles. By the remark above, we may assume that $\dict(S)=2$. Suppose first that some $Q_i$ contains a cycle $C$ and that $\dict(D)=2$. Let $\phi$ be an acyclic 2-dicolouring of $D$. Then $C$ contains an arc $uv$ such that $\phi(u)=1$ and $\phi(v)=2$.
      Since $S$ is strong and hence vertex pancyclic by Moon's theorem \cite[Theorem 1.5.1]{bang2009}, there are indices $j,k$ distinct from $i$ such that $D$ contains the subdigraph $\vec{C}_3[Q_i,Q_j,Q_k]$. If both of $Q_j,Q_k$ contain a vertex coloured $a\in\{1,2\}$ by $\phi$ then $D$ contains a monochromatic triangle so we may assume that $\phi$ uses one colour $a$ on all vertices of $Q_j$ and the other colour $b\neq a$ on all vertices of $Q_k$. Combining this with the fact that $\phi$ uses both colours of $C$ we obtain a contradiction to $\phi$ being an acyclic 2-dicolouring of $D$. 
      Assume next that each $Q_i$ is acyclic. Fix one vertex $s_i$ in each $Q_i$ and recall that these vertices induce a copy of $S$ in $D$. Now let $\psi$ be an acyclic 2-dicolouring of (this copy of) $S$. Then we can extend $\psi$ to an acyclic 2-dicolouring of $D$ by giving all vertices of $Q_i$ the colour $\psi(s_i)$.
  \end{proof}

  Notice that above we can extend any acyclic 2-dicolouring of $S$ to all of $D$ when each $Q_i$ is acyclic so we obtain the following corollary of  Theorem~\ref{thm:polytime_tour_2colours}.

  \begin{corollary}
      There is a polynomial-time algorithm for checking whether a quasi-transitive digraph $D$ has $\dict(D)\leq 2$.
  \end{corollary}

The concept of acyclic dichromatic number only makes sense for oriented graphs but, due to Theorem~\ref{thm:polytime_tour_2colours}, the following problem, which in the language of \cite{bangJGT87} can be seen as an orientation completion problem, may be of interest. Here we think of a mixed graph $M=(V,E\cup{}A)$ as corresponding to the digraph $D=(V,A^*\cup A)$ where $A^*$ is obtained from $E$ by replacing every edge in $E$ by a directed 2-cycle. Equivalently, we are looking for a spanning oriented subgraph of a digraph $D$ that we obtain by deleting precisely one arc from each 2-cycle of $D$.

\begin{problem}
    Is there a polynomial-time algorithm for deciding whether a semicomplete digraph $D=(V,A)$ contains a spanning tournament $T=(V,A')$ with $\dict(T)\leq 2$?
\end{problem}

It follows from Corollary~\ref{cor:NPCbipartite} that it is NP-complete to decide whether an oriented bipartite graph has acyclic dichromatic number 2. Hence, in general, it is NP-complete to decide whether a bipartite digraph contains a spanning oriented subdigraph with acyclic chromatic number at most 2. However, if the oriented part of a bipartite digraph $D$ is acyclic, or every arc is in a 2-cycle, then it is easy to see that $D$ always has a spanning oriented subdigraph with acyclic dichromatic number at most 2. This motivates the following problem.

\begin{problem}
    Classify the class of bipartite digraphs for which it is polynomial to decide the existence of a spanning oriented subdigraph with acyclic dichromatic number at most 2.
\end{problem}

\bibliographystyle{abbrv}
\bibliography{refs}

\end{document}